# ACCUMULATED PREDICTION ERRORS, INFORMATION CRITERIA AND OPTIMAL FORECASTING FOR AUTOREGRESSIVE TIME SERIES

By Ching-Kang Ing

*Academia Sinica and National Taiwan University*

The predictive capability of a modification of Rissanen's accumulated prediction error (APE) criterion, $\text{APE}_{\delta_n}$, is investigated in infinite-order autoregressive $(\text{AR}(\infty))$ models. Instead of accumulating squares of sequential prediction errors from the beginning, $\text{APE}_{\delta_n}$ is obtained by summing these squared errors from stage $n\delta_n$, where $n$ is the sample size and $1/n \leq \delta_n \leq 1 - (1/n)$ may depend on $n$. Under certain regularity conditions, an asymptotic expression is derived for the mean-squared prediction error (MSPE) of an AR predictor with order determined by $\text{APE}_{\delta_n}$. This expression shows that the prediction performance of $\text{APE}_{\delta_n}$ can vary dramatically depending on the choice of $\delta_n$. Another interesting finding is that when $\delta_n$ approaches 1 at a certain rate, $\text{APE}_{\delta_n}$ can achieve asymptotic efficiency in most practical situations. An asymptotic equivalence between $\text{APE}_{\delta_n}$ and an information criterion with a suitable penalty term is also established from the MSPE point of view. This offers new perspectives for understanding the information and prediction-based model selection criteria. Finally, we provide the first asymptotic efficiency result for the case when the underlying $\text{AR}(\infty)$ model is allowed to degenerate to a finite autoregression.

**1. Introduction.** In the past two decades, investigations on the accumulated prediction error (APE) [21] and its variations have attracted considerable attention among researchers from various disciplines. Prior to the early 1990s, a large number of studies focused on its consistency in selecting regression or time series models (e.g., [6, 8, 26, 27, 29]). However, since proving consistency requires assuming that the true model is included among the family of candidate models (which is rather difficult to justify in practice), recent research has focused more on understanding its statistical properties under possible model misspecification (e.g., [3, 15, 17, 20, 29, 30], among others). While a much deeper understanding of APE in cases of a misspecified model has been gained from these recent efforts, APE's prediction performance after model selection still remains unclear. This motivated the present study.

To select a model for the realization of a stationary time series, it is common to assume that the realization comes from an autoregressive moving average (ARMA) process whose AR and MA orders are known to lie within prescribed finite intervals. Then a model selection procedure is used to select orders within these intervals and thereby determine a model for the data. However, as pointed out by Shibata [25], Goldenshluger and Zeevi [5] and Ing and Wei [14], this assumption can rarely be justified in practice, and the less stringent assumption is that the time series data are observations from a linear stationary process. Following this idea, it is assumed in the sequel that observations $x_1, \ldots, x_n$ are generated by an AR($\infty$) process $\{x_t\}$, where

$$(1.1) \qquad x_t + \sum_{i=1}^{\infty} a_i x_{t-i} = e_t, \qquad t = 0, \pm 1, \pm 2, \ldots,$$

with the characteristic polynomial $A(z) = 1 + \sum_{i=1}^{\infty} a_i z^i \neq 0$ for all $|z| \leq 1$ and $\{e_t\}$ being a sequence of independent random noise variables satisfying $E(e_t) = 0$ and $E(e_t^2) = \sigma^2$ for all $t$. To predict future observations, we consider a family of approximation models $\{\mathrm{AR}(1), \ldots, \mathrm{AR}(K_n)\}$, where the maximal order $K_n$ is allowed to tend to $\infty$ as $n$ does in order to reduce approximation errors. In this framework, the APE value of model AR($k$), $1 \leq k \leq K_n$, is given by

$$(1.2) \qquad \mathrm{APE}(k) = \sum_{i=m}^{n-1} (x_{i+1} - \hat{x}_{i+1}(k))^2,$$

where $\hat{x}_{i+1}(k) = -\mathbf{x}'_i(k)\hat{\mathbf{a}}_i(k)$, $\mathbf{x}_i(k) = (x_i, \ldots, x_{i-k+1})'$, $\hat{\mathbf{a}}_i(k)$ satisfies

$$(1.3) \qquad -\hat{R}_i(k)\hat{\mathbf{a}}_i(k) = \frac{1}{i - K_n} \sum_{j=K_n}^{i-1} \mathbf{x}_j(k) x_{j+1},$$

with

$$(1.4) \qquad \hat{R}_i(k) = \frac{1}{i - K_n} \sum_{j=K_n}^{i-1} \mathbf{x}_j(k) \mathbf{x}'_j(k),$$

and $m \geq K_n + 1$ is the first integer $j$ such that $\hat{\mathbf{a}}_j(K_n)$ is uniquely defined. As observed, APE($k$) measures the performance of AR($k$) when it is used for sequential predictions. Recently, a modification of APE,

$$(1.5) \qquad \mathrm{APE}_{\delta_n}(k) = \sum_{i=n\delta_n}^{n-1} (x_{i+1} - \hat{x}_{i+1}(k))^2,$$

with $1/n \leq \delta_n \leq 1 - (1/n)$ depending on $n$, has also been considered by several authors, for example, West [30], McCracken [20] and Inoue and Kilian [15]. Since $\text{APE}_{\delta_n}$ includes the original APE as a special case, this paper focuses on $\text{APE}_{\delta_n}$. As will be shown later, the performance of $\text{APE}_{\delta_n}$ can vary dramatically depending on the choice of $\delta_n$.

In view of (1.5), it is natural to predict the next observation $x_{n+1}$ using $\hat{x}_{n+1}(\hat{k}_{n,\delta_n})$, where

(1.6) $$\hat{k}_{n,\delta_n} = \underset{1 \leq k \leq K_n}{\arg\min} \, \text{APE}_{\delta_n}(k).$$

This type of prediction, targeting future values of the observed time series, is referred to as a *same-realization prediction*. On the other hand, if the process used in estimation (or model selection) and that for prediction are independent, then it is called an *independent-realization prediction* (see [2, 16, 22] and [25]). For differences between these two types of predictions in various time series models, see [10, 11, 13, 14, 18]. The prediction performance of $\text{APE}_{\delta_n}$ after order selection is assessed using the mean-squared prediction error (MSPE) $q_n(\hat{k}_{n,\delta_n})$, where, with $\hat{k}_n = \hat{k}_n(x_1, \ldots, x_n) \in \{1, 2, \ldots, K_n\}$,

(1.7) $$q_n(\hat{k}_n) = E(x_{n+1} - \hat{x}_{n+1}(\hat{k}_n))^2.$$

There are three interrelated issues addressed in this paper. The first one focuses on the asymptotic expression for $q_n(\hat{k}_{n,\delta_n})$. To deal with this problem, in Proposition 2 (see Section 2) we establish a general theory that provides sufficient conditions under which $q_n(\hat{k}_n) - \sigma^2$ can be asymptotically expressed as a sum of two terms that measure the model complexity and the goodness of fit. This result is then applied to the case $\hat{k}_n = \hat{k}_{n,\delta_n}$ with $\delta_n$ bounded away from 1; see Theorem 1 in Section 3. A series of examples is given after Theorem 1 to illustrate its implications. In particular, it is shown in Example 1 that when the AR coefficients $\{a_i\}$ decay exponentially [which includes, but is not limited to, the $\text{ARMA}(p,q)$ model with $q > 0$ as a special case] and $\delta_n$ satisfies $\log \delta_n^{-1} = o(\log n)$, $\text{APE}_{\delta_n}$ is asymptotically efficient in the sense of (2.3). However, if the $\{a_i\}$ decay algebraically, Example 3 points out that $\text{APE}_{\delta_n}$ is no longer asymptotically efficient if $\delta_n$ is bounded away from 1. To alleviate this difficulty, Theorem 2 (also in Section 3) allows $\delta_n$ to converge at a certain rate to 1 and offers a theoretical justification for the proposed modification. In light of this result, a class of $\text{APE}_{\delta_n}$ criteria that can achieve asymptotic efficiency in both exponential and algebraic-decay cases is given; see Examples 4 and 5 after Theorem 2.

The second issue concerns the performance of the information criterion and its relation to $\text{APE}_{\delta_n}$ from the same-realization prediction point of view. The value of the information criterion for model $\text{AR}(k)$ is defined by

(1.8) $$\text{IC}_{P_n}(k) = \log \hat{\sigma}_n^2(k) + \frac{P_n k}{n},$$

where $P_n > 1$ is a positive number (possibly) depending on $n$,

$$(1.9) \qquad \hat{\sigma}_n^2(k) = \frac{1}{N} \sum_{t=K_n}^{n-1} (x_{t+1} + \hat{\mathbf{a}}_n'(k)\mathbf{x}_t(k))^2,$$

and $N = n - K_n$. Note that the AIC [1], BIC [23] and HQ criteria [7] correspond to $\mathrm{IC}_{P_n}$ with $P_n = 2, \log n$, and $c \log_2 n$, respectively, where $c > 2$ and $\log_2 n = \log(\log n)$. Equation (1.8) is referred to as an AIC-like criterion if $P_n$ is independent of $n$, and as a BIC-like criterion if $P_n \to \infty$ and $P_n = o(n)$. With the help of Proposition 2, Theorem 3 (see Section 4) gives an asymptotic expression for $q_n(\hat{k}_{n,P_n})$, where

$$(1.10) \qquad \hat{k}_{n,P_n} = \underset{1 \le k \le K_n}{\arg \min} \, \mathrm{IC}_{P_n}(k).$$

This result extends Corollary 1 of [14], which only focuses on the MSPE of the AIC-like criteria. An interesting implication of Theorem 3 is that the HQ criterion is asymptotically efficient in the exponential-decay case whereas BIC is not; see Examples 6 and 7 in Section 4. While both HQ and BIC are known to be consistent in finite-order AR models [7], these examples show that their prediction performance can differ remarkably in the AR($\infty$) case. Based on Theorems 1–3, an asymptotic equivalence between $\mathrm{IC}_{P_n}$ and $\mathrm{APE}_{\delta_n}$, with $\delta_n$ and $P_n$ satisfying (4.8), is given at the end of Section 4; see (4.9).

The third issue in which we are interested is a long-standing unresolved problem concerning time series model selection. Under the assumption that (1.1) does not degenerate to an AR model of finite order, Ing and Wei [14] recently showed that AIC, satisfying (2.3), is asymptotically efficient for same-realization predictions. However, if the order of the underlying AR model is finite, then, as mentioned previously, the BIC-like criteria (e.g., HQ and BIC) are consistent, but AIC, which asymptotically will choose an overparameterized model with positive probability, does not possess this property [24]. When the $\mathrm{APE}_{\delta_n}$ criteria are used instead, the choice between $\delta_n \to 0$ and $\delta_n \to 1$ also leads to the same difficulty; see Remark 5 in Section 3. To tackle these dilemmas, in Section 5 we first concentrate on an important special case where $\{a_i\}$ either decay exponentially or are zero for all but a finite number of $i$. It is shown in Theorem 5 that $\mathrm{IC}_{P_n}(k)$, with $P_n \to \infty$ and $P_n = o(\log n)$ and $\mathrm{APE}_{\delta_n}(k)$, with $\delta_n^{-1} \to \infty$ and $\log \delta_n^{-1} = o(\log n)$, can simultaneously achieve asymptotic efficiency over these two types of AR processes. However, if the case where $\{a_i\}$ decay algebraically is also included, then the criteria proposed by Theorem 5 fail to preserve the same optimality. A two-stage procedure, (5.1), which is a hybrid between AIC and a BIC-like criterion, is provided as a remedy. Its validity is justified theoretically in Theorem 6 (also in Section 5). Note that the results mentioned above are verified under the assumption that all moments of $e_t$ are finite [see (K.3) in Section 2]. It is made to obtain a uniform moment bound for $\hat{R}_i^{-1}(k)$ [see (B.6) in Appendix B], based on the recent work of Ing and Wei ([11], Theorem 2; [14], Proposition 1). In fact, (K.3) can be slightly relaxed at the cost of reducing the number of candidate models. However, the details are not pursued here in order to simplify the discussion. Simulation results illustrating finite sample performance of the aforementioned criteria are given in Section 6. For ease of reading, the proofs of the results in Sections 2–5 are deferred to Appendices A–D, respectively.

**2. Preliminary results.** We first list a set of assumptions that are used throughout the paper.

(K.1) Let $\{x_t\}$ be a linear process satisfying (1.1) with $A(z) = 1 + a_1 z + a_2 z^2 + \cdots \neq 0$ for $|z| \leq 1$. Furthermore, let the coefficients $\{a_i\}$ obey $\sum_{i=1}^{\infty} |i^{1/2} a_i| < \infty$.

(K.2) Let the distribution function of $e_t$ be denoted by $F_t$. There are two arbitrarily small positive numbers, $\alpha$ and $\delta_0^*$, and one arbitrarily large positive number, $C_0$, such that for all $t = \ldots, -1, 0, 1, \ldots$ and $|x - y| < \delta_0^*$,

$$|F_t(x) - F_t(y)| \leq C_0 |x - y|^\alpha.$$

(K.3) $\sup_{-\infty < t < \infty} E|e_t|^s < \infty$, $s = 1, 2, \ldots$.

(K.4) The maximal order $K_n$ satisfies

$$C_l \leq \frac{K_n^{2+\delta_1^*}}{n} \leq C_u,$$

where $\delta_1^*$, $C_l$ and $C_u$ are some prescribed positive numbers.

(K.5) $a_n \neq 0$ for infinitely many $n$.

First note that the MSPE of $\hat{x}_{n+1}(k)$, $q_n(k)$ [see (1.7)], can be expressed as

(2.1) $$\sigma^2 + E(\mathbf{f}(k) + \mathcal{S}(k))^2,$$

where

$$\mathbf{f}(k) = \mathbf{x}_n'(k) \hat{R}_n^{-1}(k) \frac{1}{N} \sum_{j=K_n}^{n-1} \mathbf{x}_j(k) e_{j+1,k}, \qquad e_{j+1,k} = x_{j+1} + \sum_{l=1}^{k} a_l(k) x_{j+1-l},$$

$$(a_1(k), \ldots, a_k(k))' = \mathbf{a}(k) = \underset{(c_1, \ldots, c_k)' \in R^k}{\arg\min} E\left(x_{k+1} + \sum_{l=1}^{k} c_l x_{k+1-l}\right)^2$$

and

$$\mathcal{S}(k) = \sum_{i=1}^{\infty} (a_i - a_i(k)) x_{n+1-i}$$

with $a_i(k) = 0$ for $i > k$. To simplify the notation, $\mathbf{a}(k)$ is sometimes viewed as an infinite-dimensional vector with undefined entries set to zero. Ing and Wei ([11], Theorem 3) obtained an asymptotic expression for $q_n(k) - \sigma^2$, which holds uniformly for all $1 \leq k \leq K_n$. This result is summarized in the following proposition.

PROPOSITION 1. *Assume that* (K.1)–(K.4) *hold. Then*

$$\lim_{n \to \infty} \max_{1 \leq k \leq K_n} \left| \frac{q_n(k) - \sigma^2}{L_n(k)} - 1 \right| = 0, \tag{2.2}$$

*where*

$$L_n(k) = \frac{k\sigma^2}{N} + \|\mathbf{a} - \mathbf{a}(k)\|_R^2,$$

*and for any infinite-dimensional vector* $\mathbf{d} = (d_1, d_2, \ldots)'$, $\|\mathbf{d}\|_R^2 = \sum_{i \leq i, j \leq \infty} d_i d_j \times \gamma_{i-j}$, *with* $\gamma_{i-j} = E(x_i x_j)$. *Also note that* $\|\mathbf{a} - \mathbf{a}(k)\|_R^2 = E(\mathcal{S}^2(k))$.

The first term of $L_n(k)$, $k\sigma^2/N$, which is proportional to $k$, can be viewed as a measure of model complexity. The second term of $L_n(k)$, $\|\mathbf{a} - \mathbf{a}(k)\|_R^2$, which decreases as $k$ increases, measures the goodness of fit. If one attempts to find an order $k$ whose corresponding predictor, $\hat{x}_{n+1}(k)$, has the minimal MSPE, then some data-driven order selection criteria are needed. An order selection criterion, $\hat{k}_n$, is said to be asymptotically efficient if $\hat{x}_{n+1}(\hat{k}_n)$ satisfies

$$\limsup_{n \to \infty} \frac{q_n(\hat{k}_n) - \sigma^2}{\min_{1 \leq k \leq K_n} q_n(k) - \sigma^2} \leq 1, \tag{2.3}$$

where $1 \leq \hat{k}_n \leq K_n$. Inequality (2.3) says that the (second-order) MSPE of the predictor with order determined by an asymptotically efficient criterion is ultimately not greater than that of the best predictor among $\{\hat{x}_{n+1}(1), \ldots, \hat{x}_{n+1}(K_n)\}$. In view of (2.2), (2.3) is equivalent to

$$\limsup_{n \to \infty} \frac{q_n(\hat{k}_n) - \sigma^2}{L_n(k_n^*)} \leq 1, \tag{2.4}$$

where $L_n(k_n^*) = \min_{1 \leq k \leq K_n} L_n(k)$.

Let $\mathrm{OS}_n(k)$ be an order selection function and

$$\hat{k}_{n,\mathrm{OS}} = \arg\min_{1 \leq k \leq K_n} \mathrm{OS}_n(k) \tag{2.5}$$

be the selected order. We shall provide sufficient conditions under which $q_n(\hat{k}_{n,\mathrm{OS}}) - \sigma^2$ can be asymptotically expressed in terms of the $L_n(\cdot)$ function. Define

$$L_{n,D_n}(k) = \frac{(D_n - 1)k\sigma^2}{N} + \|\mathbf{a} - \mathbf{a}(k)\|_R^2, \tag{2.6}$$

where $D_n > 1$, and

(2.7) $$k^*_{n,D_n} = \underset{1 \leq k \leq K_n}{\arg\min}\, L_{n,D_n}(k).$$

PROPOSITION 2. *Assume that* (K.1)–(K.4) *hold. If there exists a sequence of positive numbers* $\{D_n\}$, *with* $\liminf_{n \to \infty} D_n > 1$, *such that*

(2.8) $$\lim_{n \to \infty} (D_n - 1) \frac{E(\mathcal{S}(\hat{k}_{n,\mathrm{OS}}) - \mathcal{S}(k^*_{n,D_n}))^2}{L_{n,D_n}(k^*_{n,D_n})} = 0$$

*and*

(2.9) $$\lim_{n \to \infty} (D_n - 1) \frac{E(\mathbf{f}(\hat{k}_{n,\mathrm{OS}}) - \mathbf{f}(k^*_{n,D_n}))^2}{L_{n,D_n}(k^*_{n,D_n})} = 0,$$

*where* $\mathcal{S}(k)$ *and* $\mathbf{f}(k)$ *are defined after* (2.1), *then*

(2.10) $$\lim_{n \to \infty} \frac{q_n(\hat{k}_{n,\mathrm{OS}}) - \sigma^2}{L_n(k^*_{n,D_n})} = 1.$$

*Moreover, if*

(2.11) $$\lim_{n \to \infty} \frac{L_n(k^*_{n,D_n})}{L_n(k^*_n)} = 1,$$

*then* (2.3) [(2.4)] *holds for* $\hat{k}_n = \hat{k}_{n,\mathrm{OS}}$.

REMARK 1. If (K.1), (K.5), $\sup_{-\infty < t < \infty} E|e_t|^4 < \infty$ and $K_n = o(n^{1/2})$ are assumed, and (2.8) and (2.9) are replaced with

(2.12) $$\underset{n \to \infty}{\text{p-lim}}(D_n - 1) \frac{L_n(\hat{k}_{n,\mathrm{OS}}) - L_n(k^*_{n,D_n})}{L_{n,D_n}(k^*_{n,D_n})} = 0,$$

then it is shown in Appendix A that

(2.13) $$\underset{n \to \infty}{\text{p-lim}} \frac{E\{(y_{n+1} - \hat{y}_{n+1}(\hat{k}_{n,\mathrm{OS}}))^2 | x_1, \ldots, x_n\} - \sigma^2}{L_n(k^*_{n,D_n})} = 1,$$

where $y_{n+1}$ is the future value of $\{y_1, \ldots, y_n\}$, which is a realization from an independent copy of $\{x_t\}$, and $\hat{y}_{n+1}(k) = -\mathbf{y}'_n(k)\hat{a}_n(k)$ with $\mathbf{y}'_n(k) = (y_n, \ldots, y_{n+1-k})$. Note that (2.13) gives an asymptotic expression for the (conditional) MSPE of $\hat{k}_{n,\mathrm{OS}}$ in independent-realization settings. For further discussion, see Remark 6 in Section 6.

Proposition 2 asserts that if $\hat{k}_{n,\mathrm{OS}}$ is sufficiently close to $k^*_{n,D_n}$ in the sense of (2.8) and (2.9), then $q_n(\hat{k}_{n,\mathrm{OS}}) - \sigma^2$ has the asymptotic expression $L_n(k^*_{n,D_n})$. In addition, if (2.11) also holds, then $\hat{k}_{n,\text{OS}}$ is asymptotically efficient. Proposition 2 plays a prominent role in justifying $\text{APE}_{\delta_n}$'s and $\text{IC}_{P_n}$'s asymptotic (in)efficiency in various situations; see Sections 3–5. To apply Proposition 2, it is important to determine a penalty term $D_n$ associated with the selection criterion $\text{OS}_n(k)$ and then justify (2.8) and (2.9) under suitable assumptions. For the first task, it is shown in Section 3 that the $D_n$ associated with $\text{APE}_{\delta_n}(k)$ is

$$(2.14) \qquad D_{\text{APE}_{\delta_n}} = 1 + \frac{N \log \delta_n^{-1}}{n(1-\delta_n)}.$$

According to Appendix C, the $D_n$ associated with $\text{IC}_{P_n}(k)$ and

$$(2.15) \qquad \frac{S_n^{(P_n)}(k)}{N} = \left(1 + \frac{P_n k}{N}\right)\hat{\sigma}_n^2(k)$$

is $P_n$. To facilitate the second task, inspired by Ing and Wei ([14], (3.9)), assumption (K.6) (see below) is frequently used in the rest of this paper. As will be seen in Appendices B–D, it is introduced to deal with the complicated dependency conditions among the selected order, estimated parameters and future observations. (Note that in the independent-realization settings, the future value to be predicted is independent of the selected order and estimated parameters.)

(K.6) For any $\xi > 0$, there are a nonnegative exponent $0 \leq \theta = \theta(\xi) < 1$ and a positive number $M = M(\xi)$ such that

$$(2.16) \qquad \liminf_{n\to\infty} R_n(\xi, \theta, M) > 0,$$

where, with $K_n$ obeying (K.4), $D_n$ satisfying $\liminf D_n > 1$ and $D_n = o(n)$, and $A_{D_n,\theta,M} = \{k : 1 \leq k \leq K_n, |k - k^*_{n,D_n}| \geq M(k^*_{n,D_n})^\theta\}$,

$$R_n(\xi, \theta, M) = \min_{k \in A_{D_n,\theta,M}} (k^*_{n,D_n})^\xi \frac{N\{L_{n,D_n}(k) - L_{n,D_n}(k^*_{n,D_n})\}}{(D_n - 1)|k - k^*_{n,D_n}|}.$$

If $\{x_t\}$ is an AR process of finite order [viz., (K.5) does not hold], then (2.16) automatically holds. On the other hand, if (K.5) holds instead, by arguments similar to those in Examples 1 and 2 and the Appendix of Ing and Wei [14], it can also be shown that (2.16) is satisfied in the following cases: (a) the exponential-decay case,

$$(2.17) \qquad C_1 k^{-\theta_1} e^{-\beta k} \leq \|\mathbf{a} - \mathbf{a}(k)\|_R^2 \leq C_2 k^{\theta_1} e^{-\beta k},$$

where $C_2 \geq C_1 > 0, \theta_1 \geq 0$ and $\beta > 0$ [note that if (K.1) is assumed, then (2.17) is equivalent to $C_1^* k^{-\theta_1} e^{-\beta k} \leq \sum_{i \geq k} a_i^2 \leq C_2^* k^{\theta_1} e^{-\beta k}$, for some $C_2^* \geq C_1^* > 0$]; and (b) the algebraic-decay case,

$$(2.18) \qquad (C_3 - M_1 k^{-\xi_1}) k^{-\beta} \leq \|\mathbf{a} - \mathbf{a}(k)\|_R^2 \leq (C_3 + M_1 k^{-\xi_1}) k^{-\beta},$$

where $C_3, M_1 > 0$, $\xi_1 \geq 2$ and $\beta > 1 + \delta_1^*$ [recall that $\delta_1^*$ is defined in (K.4)]. These facts reveal that (2.16) is quite reasonable from both practical and theoretical points of view, since it includes the ARMA model (which is the most used short-memory time series model by far) and the AR($\infty$) model with algebraically decaying coefficients (which is of much theoretical interest in the context of model selection) as special cases.

It is worth mentioning that when (K.1)–(K.6) are assumed, (2.8) and (2.9) were verified by Ing and Wei ([14], (5.75) and (5.74), resp.) in the special case where $\hat{k}_{n,\text{OS}} = \arg\min_{1 \leq k \leq K_n} S_n^{(2)}(k)$ and $D_n = 2$. Using similar arguments and assumptions, it can be shown that (2.8) and (2.9) are still valid for $\hat{k}_{n,\text{OS}} = \hat{k}_{n,P_n}$ [defined in (1.10)], $D_n = P_n$ and $1 < P_n = \alpha < \infty$ independent of $n$. However, to justify (2.8) and (2.9) in the case $\hat{k}_{n,\text{OS}} = \hat{k}_{n,\delta_n}$ and $D_n = D_{\text{APE}_{\delta_n}}$ or in the case $\hat{k}_{n,\text{OS}} = \hat{k}_{n,P_n}$ and $D_n = P_n \to \infty$, a much more delicate analysis is required. This problem is tackled in the next two sections. Note that in the finite-order AR models, (2.8) and (2.9) can also be verified for these two types of criteria under suitable assumptions; see Section 5 for more details.

As observed in Proposition 2, (2.11) is an important key to the asymptotic efficiency. It holds if the penalty term $D_n$ satisfies

$$\lim_{n \to \infty} D_n = 2. \tag{2.19}$$

To see this, first observe that

$$\max_{1 \leq k \leq K_n} \left| \frac{L_{n,D_n}(k)}{L_n(k)} - 1 \right| \to 0 \tag{2.20}$$

if $D_n \to 2$. The result (2.20) and the fact that

$$1 \leq \frac{L_n(k_{n,D_n}^*)}{L_n(k_n^*)} = \frac{L_n(k_{n,D_n}^*)/L_{n,D_n}(k_{n,D_n}^*)}{L_n(k_n^*)/L_{n,D_n}(k_n^*)} \frac{L_{n,D_n}(k_{n,D_n}^*)}{L_{n,D_n}(k_n^*)}$$

$$\leq \frac{L_n(k_{n,D_n}^*)/L_{n,D_n}(k_{n,D_n}^*)}{L_n(k_n^*)/L_{n,D_n}(k_n^*)} \tag{2.21}$$

yield (2.11). When $\|\mathbf{a} - \mathbf{a}(k)\|_R^2$ decays exponentially or (K.5) is violated, (2.11) can hold without (2.19); see Example 1 in Section 3 and the proof of Theorem 4 in Appendix D. For some other interesting discussion regarding (2.11), see Examples 2 and 3 in Section 3 and Examples 6–8 in Section 4.

**3. The MSPE of APE$_{\delta_n}$ in AR($\infty$) processes.** This section provides asymptotic expressions for $q_n(\hat{k}_{n,\delta_n}) - \sigma^2$. Without loss of generality, $n\delta_n, 1/n \leq \delta_n \leq 1 - (1/n)$, is assumed to be a positive integer. First note that

$$\text{APE}_{\delta_n}(k) = \sum_{i=n\delta_n}^{n-1} (x_{i+1} + \mathbf{x}_i'(k)\hat{\mathbf{a}}_i(k))^2$$

(3.1)
$$= \sum_{i=n\delta_n}^{n-1} \{e_{i+1} + \hat{e}_{i,k} + (e_{i+1,k} - e_{i+1})\}^2,$$

where $\hat{e}_{i,k} = \mathbf{x}'_i(k)(\hat{\mathbf{a}}_i(k) - \mathbf{a}(k))$ and $e_{i+1,k}$ is defined after (2.1). Following Lai and Wei ([19], (2.7)),

(3.2)
$$\sum_{i=n\delta_n}^{n-1} \hat{e}_{i,k}^2 = \sum_{i=n\delta_n}^{n-1} h_i(k)e_{i+1,k}^2 + Q_{n\delta_n}(k) - Q_n(k) + \sum_{i=n\delta_n}^{n-1} h_i(k)\hat{e}_{i,k}^2$$
$$- 2\sum_{i=n\delta_n}^{n-1} (1 - h_i(k))\hat{e}_{i,k}e_{i+1,k},$$

where
$$h_i(k) = \mathbf{x}'_i(k)\left(\sum_{j=K_n}^{i} \mathbf{x}_j(k)\mathbf{x}'_j(k)\right)^{-1}\mathbf{x}_i(k)$$

and
$$Q_i(k) = \left(\sum_{j=K_n}^{i-1} \mathbf{x}_j(k)e_{j+1,k}\right)'\left(\sum_{j=K_n}^{i-1} \mathbf{x}_j(k)\mathbf{x}'_j(k)\right)^{-1}\left(\sum_{j=K_n}^{i-1} \mathbf{x}_j(k)e_{j+1,k}\right).$$

On substituting (3.2) into (3.1), one obtains

(3.3)
$$\text{APE}_{\delta_n}(k) = \sum_{i=n\delta_n}^{n-1} e_{i+1}^2 + \left\{\sum_{i=n\delta_n}^{n-1} h_i(k)e_{i+1,k}^2 - k\sigma^2 \log \delta_n^{-1}\right\}$$
$$+ Q_{n\delta_n}(k) - Q_n(k) + \sum_{i=n\delta_n}^{n-1} h_i(k)\hat{e}_{i,k}^2 + 2\sum_{i=n\delta_n}^{n-1} h_i(k)\hat{e}_{i,k}e_{i+1,k}$$
$$+ \sum_{i=n\delta_n}^{n-1} \{(e_{i+1,k} - e_{i+1})^2 - \|\mathbf{a} - \mathbf{a}(k)\|_R^2\}$$
$$+ 2\sum_{i=n\delta_n}^{n-1} (e_{i+1,k} - e_{i+1})e_{i+1} + n(1 - \delta_n)L_{n,D_{\text{APE}_{\delta_n}}}(k),$$

where $D_{\text{APE}_{\delta_n}}$ is defined in (2.14). When $\delta_n$ is bounded away from 1, Theorem 1 below provides sufficient conditions under which (2.8) and (2.9) hold for $\hat{k}_{n,\text{OS}} = \hat{k}_{n,\delta_n}$ and $D_n = D_{\text{APE}_{\delta_n}}$. As a result, an asymptotic expression for $q_n(\hat{k}_{n,\delta_n}) - \sigma^2$ is obtained. Note that the relation $D_n = D_{\text{APE}_{\delta_n}}$ is used in the rest of this section.

THEOREM 1. *Assume that* (K.1)–(K.6) *hold and* $1/n \leq \delta_n \leq 1 - (1/n)$ *satisfies*

$$\limsup_{n \to \infty} \delta_n < 1 \tag{3.4}$$

*and*

$$0 < \liminf_{n \to \infty} n^{\theta_3} \delta_n \leq \infty, \tag{3.5}$$

*where* $0 < \theta_3 < \delta_1^*/(2+\delta_1^*)$. *Moreover, if for some* $0 \leq \theta = \theta(\xi) < 1$ *and* $\eta > 0$,

$$\lim_{n \to \infty} \frac{D_n - 1}{(k_{n,D_n}^*)^{1-\theta(1+\eta)}} = 0, \tag{3.6}$$

*where* $\theta$ *is obtained from* (K.6) *when*

$$0 < \xi < \min\{1/2, \{(2+\delta_1^*)(1-\theta_3)/2\} - 1\}, \tag{3.7}$$

*then* (2.8) *and* (2.9) *are true for* $\hat{k}_{n,OS} = \hat{k}_{n,\delta_n}$. *Hence,*

$$\lim_{n \to \infty} \frac{q_n(\hat{k}_{n,\delta_n}) - \sigma^2}{L_n(k_{n,D_n}^*)} = 1. \tag{3.8}$$

REMARK 2. Note that $\theta$ in (3.6) is not uniquely determined. In order for (3.6) to be less stringent, $\theta$ can be chosen as small as possible; see Examples 1 and 3 below for more details.

REMARK 3. Consider the following assumption:

(K.6′) For any $\xi > 0$, there is a subsequence of $\{n\}$, $\{n_l\}$, such that

$$\liminf_{l \to \infty} R_{n_l}(\xi, 0, 1) > 0. \tag{3.9}$$

[Recall that $R_n(\xi, \theta, M)$ is defined in (K.6).] If, in place of (K.6), (K.6′) is assumed in Theorem 1, then it can be shown that (3.8) remains valid for $n = n_l$, without imposing (3.6). This finding is applied in Example 2 below to illustrate that the $\text{APE}_{\delta_n}$ criteria, with $\delta_n$ decreasing to 0 at a polynomial rate, perform poorly in the case where the AR coefficients decay exponentially fast.

REMARK 4. Under (K.5), it is not difficult to see that $k_{n,D_n}^* \to \infty$ as $n \to \infty$. Therefore, when $0 < \delta_n = \delta < 1$ is fixed with $n$, (3.6) automatically holds.

The following examples help gain a better understanding of Theorem 1.

EXAMPLE 1. Assume that (K.1)–(K.4) hold and the AR coefficients satisfy

$$C_1 e^{-\beta k} \leq \|\mathbf{a} - \mathbf{a}(k)\|_R^2 \leq C_2 e^{-\beta k}, \tag{3.10}$$

where $0 < C_1 \leq C_2 < \infty$ and $\beta > 0$. Note that (3.10) is satisfied by any causal and invertible ARMA$(p,q)$ model with $q > 0$. As mentioned in Section 2, when (K.1) is assumed, (3.10) can also be expressed as $C_1' e^{-\beta k} \leq \sum_{i \geq k} a_i^2 \leq C_2' e^{-\beta k}$ for some $0 < C_1' \leq C_2' < \infty$. We shall show in this example that (3.8) follows if $\delta_n$ satisfies (3.4) and

$$\log \delta_n^{-1} = o(\log n). \tag{3.11}$$

Condition (3.11) guarantees (3.5). It can be shown to be equivalent to $n^\nu \delta_n \to \infty$ for all $\nu > 0$. Obviously, (3.4) and (3.11) are satisfied if $0 < \delta_n = \delta < 1$ is independent of $n$ or if $\delta_n^{-1} = (\log n)^{\nu_1}$ for some $\nu_1 > 0$. Therefore, in view of Theorem 1 and the discussion given after (K.6), it remains to verify (3.6). By (3.10) and an argument similar to that used in the Appendix of [14], for some $C_1 > 0$,

$$\frac{1}{\beta} \log n - \frac{1}{\beta} \log(D_n - 1) - C_1 \leq k_{n,D_n}^* \tag{3.12}$$
$$\leq \frac{1}{\beta} \log n - \frac{1}{\beta} \log(D_n - 1) + C_1,$$

and for any $\xi > 0$, (2.16) holds for $\theta = 0$ and some $M > 0$ (or for any $0 < \theta < 1$ and any $M > 0$). As a result, (3.6) holds for $\theta = 0$ and $\eta > 0$, and hence (3.8) follows. Moreover, (3.10) and the same argument used to prove (A.1) of [14] yield that for some $C_2 > 0$,

$$\frac{1}{\beta} \log n - C_2 \leq k_n^* \leq \frac{1}{\beta} \log n + C_2. \tag{3.13}$$

According to (3.10)–(3.13), we obtain (2.11), which together with (3.8) implies that $\text{APE}_{\delta_n}$, with $\delta_n$ satisfying (3.4) and (3.11), is asymptotically efficient.

EXAMPLE 2. This example is given to indicate that if $\delta_n$ decays to 0 at a polynomial rate, then $\text{APE}_{\delta_n}$ cannot be asymptotically efficient even in the exponential-decay case. More specifically, assume that (K.1)–(K.4) are satisfied,

$$\delta_n = C_1 n^{-\theta_3}, \tag{3.14}$$

where $C_1 > 0$ and $0 < \theta_3 < \delta_1^*/(2+\delta_1^*)$, and the AR coefficients obey a special case of (3.10),

$$\|\mathbf{a} - \mathbf{a}(k)\|_R^2 = C_2 e^{-\beta k}(1 + G_k), \tag{3.15}$$

where $|G_k| < 1$, $G_k \to 0$ as $k \to \infty$ and $C_2$ and $\beta$ are some positive numbers. By (3.15) and an argument similar to that used in Case II of [25], page 162, for sufficiently large $n$,

$$(3.16) \qquad k_{n,D_n}^* = m_{1,n} \text{ or } m_{1,n} + 1,$$

where $m_{1,n}$ is the largest integer $\leq m_{1,n}^* = (1/\beta)\log[C_2 N\beta/((D_n-1)\sigma^2)]$. Define $x_{1,n} = m_{1,n}^* - m_{1,n}$ and $z = (1/\beta)\log[\beta/(1-e^{-\beta})]$. Since $0 < z < 1$, there is a positive number $\kappa$ such that $0 < z - \kappa < z + \kappa < 1$. Define a set of positive integers $A_\kappa = \{n : |x_{1,n} - z| > \kappa, n = 1, 2, \ldots\}$. Then, it can be shown that $A_\kappa$ contains infinitely many elements. Moreover, for any $\xi > 0$ and any sequence of positive integers $\{n_l\} \subseteq A_\kappa$, (3.9) holds. Therefore, according to Remark 3, (3.8) is valid for $n = n_l$. By analogy with (3.16), for sufficiently large $n$,

$$(3.17) \qquad k_n^* = m_{2,n} \text{ or } m_{2,n} + 1,$$

where $m_{2,n}$ is the largest integer $\leq (1/\beta)\log(C_2 N\beta/\sigma^2)$. (3.14)–(3.17) yield

$$(3.18) \qquad \lim_{n\to\infty} \frac{L_n(k_{n,D_n}^*)}{L_n(k_n^*)} > 1,$$

which, together with (3.8) (with $n = n_l$) gives

$$(3.19) \qquad \limsup_{n\to\infty} \frac{q_n(\hat{k}_{n,\delta_n}) - \sigma^2}{L_n(k_n^*)} > 1.$$

As a result, $\text{APE}_{\delta_n}$, with $\delta_n$ satisfying (3.14), fails to achieve (2.4) in the exponential-decay case.

EXAMPLE 3. This example investigates the prediction performance of $\text{APE}_{\delta_n}$ in the algebraic-decay case (2.18). If (2.18), (3.4) and (3.5) are satisfied, then the same argument as the one in Example 2 of [14] yields that

$$(3.20) \qquad k_{n,D_n}^* = (NC_3\beta(D_n-1)^{-1}\sigma^{-2})^{1/(\beta+1)} + O(1),$$

and for any $\xi > 0$, (2.16) holds for any $1 - \min\{\xi, 1\} < \theta < 1$ and any $M > 0$. These facts and (3.5) guarantee that (3.6) is valid for $1 - \min\{\xi, 1\} < \theta < 1$ and $0 < \eta < (1-\theta)/\theta$, where $\xi$ satisfies (3.7). Consequently, when (K.1)–(K.4), (2.18), (3.4) and (3.5) are assumed, (3.8) is ensured by Theorem 1. By (A.9) of [14],

$$(3.21) \qquad k_n^* = (NC_3\beta\sigma^{-2})^{1/(\beta+1)} + O(1).$$

This, (2.18), (3.4), (3.5) and (3.20) imply that

$$(3.22) \qquad \liminf_{n\to\infty} \frac{L_n(k_{n,D_n}^*)}{L_n(k_n^*)} > 1.$$

According to (3.8) and (3.22), the $\text{APE}_{\delta_n}$ is not asymptotically efficient in this case.

As can be seen from Example 3, due to violation of (2.11), $\text{APE}_{\delta_n}$, with $\delta_n$ bounded away from 1, is not asymptotically efficient in the algebraic-decay case. In view of (2.19)–(2.21), a natural remedy for this difficulty is to let $\delta_n \to 1$. However, problems can still occur if $\delta_n$ converges "too fast" to 1. To see this, let $\delta_n = 1 - (1/n)$. Then $\text{APE}_{\delta_n}(k) = (x_n - \hat{x}_n(k))^2$. Since models are determined only by the last period's prediction errors, it does not seem possible to establish any (asymptotically) optimal selection result in this case. To resolve this dilemma, some suitable choices of $\delta_n$ are introduced in Theorem 2. Some examples are also given after the theorem to help gain further insight into it.

THEOREM 2. *Assume that* (K.1)–(K.6) *hold and* $1/n \leq \delta_n \leq 1 - (1/n)$ *satisfies* $\lim_{n\to\infty} \delta_n = 1$. *Moreover, if either of the following conditions holds, then* (2.8) *and* (2.9) *are valid for* $\hat{k}_{n,\text{OS}} = \hat{k}_{n,\delta_n}$.

(i) $\lim_{n\to\infty} k^*_{n,D_n}/n^{\theta_3} = 0$ *for any* $\theta_3 > 0$ *and* $(1-\delta_n)^{-1} = O(k^{*\xi_2}_{n,D_n})$ *for some* $0 < \xi_2 < 1/2$.

(ii) $(1-\delta_n)^{-1} = O(k^{*\xi_2}_{n,D_n})$ *for some* $0 < \xi_2 < \min\{1/2, \delta_1^*/2\}$.

*Consequently,* $APE_{\delta_n}$ *is asymptotically efficient in the sense of* (2.3) [(2.4)].

In light of Theorem 2, the following examples demonstrate how to choose $\delta_n$ such that the resulting $\text{APE}_{\delta_n}$ is asymptotically efficient in both the exponential- and algebraic-decay cases.

EXAMPLE 4. Assume that (K.1)–(K.4) hold and the AR coefficients obey (2.17). Although Example 1 shows that when $\theta_1$ in (2.17) is equal to 0, $\text{APE}_{\delta_n}$, with $\delta_n$ satisfying (3.4) and (3.11), is asymptotically efficient, it is unclear whether this result still holds for $\theta_1 > 0$. Fortunately, this difficulty can be bypassed by letting

$$\delta_n = 1 - C_1(\log n)^{-r}, \qquad (3.23)$$

with $C_1 > 0$ and $0 < r < 1/2$. First note that under (2.17) and (3.23), the same argument as in Example 1 of [14] yields that for some $C_2 > 0$,

$$(3.24) \qquad \frac{1}{\beta}\log n - C_2 \log_2 n \leq k^*_{n,D_n} \leq \frac{1}{\beta}\log n + C_2 \log_2 n,$$

and for any $\xi > 0$, (2.16) holds for any $0 < \theta < 1$ and any $M > 0$. Moreover, since condition (i) of Theorem 2 is ensured by (3.23) and (3.24), $\text{APE}_{\delta_n}$, with $\delta_n$ satisfying (3.23), is asymptotically efficient under (2.17).

EXAMPLE 5. This example shows that if $\delta_n$ satisfies (3.23) with $C_1 > 0$ and $0 < r < \infty$, then the corresponding $\text{APE}_{\delta_n}$ is asymptotically efficient under the algebraic-decay case (2.18). To see this, first note that following the same line of reasoning as in Example 2 of [14], (3.20) is still valid, and for any $\xi > 0$, (2.16) holds for any $1 - \min\{\xi, 1\} < \theta < 1$ and any $M > 0$. In addition, since condition (ii) of Theorem 2 is ensured by (3.20) and the condition imposed on $\delta_n$, the desired result follows from Theorem 2.

Examples 4 and 5 suggest that to achieve asymptotic efficiency through $\text{APE}_{\delta_n}$ in both the exponential- and algebraic-decay cases, $\delta_n$ can be chosen to satisfy (3.23) with $C_1 > 0$ and $0 < r < 1/2$. However, the question of how to determine the best $C_1$ and $r$ seems difficult to answer from a finite sample point of view. For some simulation results illustrating $\text{APE}_{\delta_n}$'s performance in finite samples, see Section 5. We close this section with two remarks concerning the performance of $\text{APE}_{\delta_n}$ in finite-order AR models and for independent-realization predictions.

REMARK 5. When (1.1) degenerates to an $\text{AR}(p_0)$ model with $1 \leq p_0 < \infty$, it can be shown that $\hat{k}_{n,\delta_n}$, with $\liminf_{n \to \infty} \delta_n > 0$, is not a consistent estimator of $p_0$ (e.g., [15]). On the other hand, if $\delta_n \to 0$ at a certain rate, then the corresponding $\text{APE}_{\delta_n}$ is consistent and asymptotically efficient (see Appendix D). Since these results and Theorem 2 offer totally different suggestions for choosing $\delta_n$, it becomes very challenging to achieve asymptotic efficiency through $\text{APE}_{\delta_n}$ when (1.1) is allowed to degenerate to a finite autoregression. In Section 5, some selection criteria to remedy this difficulty are proposed.

REMARK 6. Note that the $\text{APE}_{\delta_n}$ described in Theorem 2 is also asymptotically efficient for independent-realization predictions. By Corollary B.1 (see Appendix B),

$$(3.25) \qquad \operatorname*{p-lim}_{n \to \infty} \frac{L_{n,D_n}(\hat{k}_{n,\delta_n})}{L_{n,D_n}(k^*_{n,D_n})} - 1 = 0.$$

Armed with (3.25) and (2.19)–(2.21), it can be shown that (2.12) holds for $\hat{k}_{n,\text{OS}} = \hat{k}_{n,\delta_n}$. Consequently, Remark 1 and (2.11) guarantee that

$$(3.26) \qquad \operatorname*{p-lim}_{n \to \infty} \frac{E\{(y_{n+1} - \hat{y}_{n+1}(\hat{k}_{n,\delta_n}))^2 | x_1, \ldots, x_n\} - \sigma^2}{L_n(k^*_n)} = 1,$$

which gives the claimed result. For more details on the definition of asymptotic efficiency in independent-realization settings, see [2, 16] and [25].

**4. The MSPE of $\text{IC}_{P_n}$ in $\text{AR}(\infty)$ processes.** In this section, prediction performance of the information criterion $\text{IC}_{P_n}(k), P_n > 1$, is investigated. When $P_n$ is independent of $n$, Ing and Wei ([14], Corollary 1) obtained an asymptotic expression for $q_n(\hat{k}_{n,P_n}) - \sigma^2$, where $\hat{k}_{n,P_n}$, defined in (1.10), is the minimizer of $\text{IC}_{P_n}(k)$, with $1 \le k \le K_n$ and $K_n$ satisfying (K.4). Theorem 3 below extends Ing and Wei's result to the case where $P_n$ is allowed to tend to $\infty$ with $n$. Note that the relation $D_n = P_n$ is used throughout this section.

THEOREM 3. *Let* (K.1)–(K.6) *hold and $P_n$ satisfy*

(4.1) $$\liminf_{n \to \infty} P_n > 1$$

*and*

(4.2) $$P_n = O(n^{\theta_3}),$$

*for some $0 < \theta_3 < (1 + \delta_1^*)/(4 + 2\delta_1^*)$. Moreover, if* (3.6) *holds with* (3.7) *replaced by*

(4.3) $\quad 0 < \xi < \min\{1/2, \delta_1^*/2, \{(1/2) - \theta_3\}(2 + \delta_1^*), (1 + \delta_1^*) - 2\theta_3(2 + \delta_1^*)\},$

*then* (2.8) *and* (2.9) *are true for $\hat{k}_{n,OS} = \hat{k}_{n,P_n}$. Consequently,*

(4.4) $$\lim_{n \to \infty} \frac{q_n(\hat{k}_{n,P_n}) - \sigma^2}{L_n(k_{n,D_n}^*)} = 1.$$

REMARK 7. If in Theorem 3, (K.6′) holds instead of (K.6), then it can be shown that (4.4) is still valid for $n = n_l$ [note that $n_l$ is defined in (K.6′)]. In this case, condition (3.6) is not required. This result can be applied to verify that BIC is not asymptotically efficient in the exponential-decay case; see Example 7 for more details.

REMARK 8. Since (K.5) is assumed, (3.6) holds automatically if $P_n = O(1)$.

The following examples illustrate implications of Theorem 3. Special emphasis is placed on comparing the predictive capabilities of three well-known information criteria, AIC, HQ and BIC, in various situations.

EXAMPLE 6. Assume that (K.1)–(K.4) hold and the AR coefficients satisfy (3.10). We shall show that $\text{IC}_{P_n}(k)$, with $P_n$ satisfying (4.1) and

(4.5) $$P_n = o(\log n),$$

is asymptotically efficient. Therefore, the AIC and HQ criteria are asymptotically efficient in this case. To see this, first note that the same reasoning as in Example 1 yields (3.12) and that for any $\xi > 0$, (2.16) holds for $\theta = 0$ and some $M > 0$ (or for any $0 < \theta < 1$ and any $M > 0$). These results and (4.5) imply that (3.6) holds for $\theta = 0$ and $\eta > 0$. According to Theorem 3, (4.4) follows. Consequently, the claimed result is ensured by (4.4) and by observing that (2.11) is valid under (3.10), (3.12), (3.13) and (4.5).

EXAMPLE 7. This example illustrates that an information criterion cannot be asymptotically efficient in the exponential-decay case when the weight for penalizing the number of regressors in the model is "too strong." To see this, let (K.1)–(K.4) and (3.15) be satisfied and

$$(4.6) \qquad P_n = C_1 (\log n)^{C_2},$$

where $C_1 > 0$ and $C_2 \geq 1$. Under these assumptions, (3.16) is obtained and (3.9) holds for any $\xi > 0$ and any sequence $\{n_l\} \subseteq A_\kappa$, where $A_\kappa$ is defined as in Example 2. By Remark 7, (4.4) is valid for $n = n_l$. Moreover, since (3.15)–(3.17) and (4.6) yield (3.18), it is concluded that $\text{IC}_{P_n}(k)$, with $P_n$ given by (4.6), is not asymptotically efficient. One important implication of this example is that BIC is not asymptotically efficient in the algebraic-decay case.

EXAMPLE 8. Consider the algebraic-decay case, (2.18). Let $P_n$ satisfy (4.1) and

$$(4.7) \qquad P_n = O((\log n)^{C_1}),$$

for some $C_1 > 0$. By an argument similar to that used in Example 3, one obtains (3.20), and for any $\xi > 0$, (2.16) holds for any $1 - \min\{\xi, 1\} < \theta < 1$ and any $M > 0$. These facts and (4.7) imply that (3.6) is valid for $1 - \min\{\xi, 1\} < \theta < 1$ and $0 < \eta < (1-\theta)/\theta$, where $\xi$ satisfies (4.3). As a result, (4.4) follows from Theorem 3. Moreover, by (2.19)–(2.21), (2.11) holds when $\lim_{n\to\infty} P_n = 2$; and by (2.18), (3.20) and (3.21), $\limsup_{n\to\infty} L_n(k^*_{n,D_n})/L_n(k^*_n) > 1$ when $\lim_{n\to\infty} P_n \neq 2$. This observation and (4.4) imply that AIC and $\text{AIC}_C$ [9] are asymptotically efficient in the algebraic-decay case, (2.18), whereas HQ, BIC and any information criterion with $\lim_{n\to\infty} P_n \neq 2$ are not.

As a final remark, note that when the conditions imposed by Theorems 1 and 3 (or Theorems 2 and 3) hold and

$$(4.8) \qquad \lim_{n\to\infty} \frac{\log \delta_n^{-1}}{(1-\delta_n)(P_n - 1)} = 1,$$

then

$$(4.9) \qquad \lim_{n\to\infty} \frac{E(x_{n+1} - \hat{x}_{n+1}(\hat{k}_{n,P_n}))^2 - \sigma^2}{E(x_{n+1} - \hat{x}_{n+1}(\hat{k}_{n,\delta_n}))^2 - \sigma^2} = 1.$$

Instead of attempting to achieve a certain asymptotic optimality for prediction, (4.8) and (4.9) are interesting in that (4.8) can be used to connect the sequence $\delta_n$ defining $\text{APE}_{\delta_n}$ with a corresponding parameter estimation penalty weight sequence, $P_n = 1 + (1-\delta_n)^{-1} \log \delta_n^{-1}$, in such a way that $\hat{x}_{n+1}(\hat{k}_{n,P_n})$ and $\hat{x}_{n+1}(\hat{k}_{n,\delta_n})$ have the same asymptotic same-realization prediction (in)efficiency, as observed in (4.9). And, conversely, a sequence $P_n$ implicitly determines a sequence $\delta_n$ through this same relation, which yields identical asymptotic (in)efficiency for $\text{IC}_{P_n}$ and $\text{APE}_{\delta_n}$. This connection not only imparts to the $\text{APE}_{\delta_n}$ criteria the deep foundations of the information criteria, but also endows the information criteria with an on-line prediction meaning. For a related result, Wei ([29], Theorem 4.2.2), under (1.1) and certain moment conditions on $e_t$ (which can be verified for the normal distribution), established an algebraic connection between BIC and APE, $\log(\text{APE}(k)/n) = \text{BIC}(k) + o(\log n/n)$ a.s., where $k$ is a positive integer and *fixed* with $n$. Therefore, except for the $o(\log n/n)$ term, the logarithm of $\text{APE}(k)/n$ is (a.s.) identical to $\text{BIC}(k)$. Hannan, McDougall and Poskitt [6] also obtained the same result in a stationary $\text{AR}(p_0)$ model with $p_0 < \infty$ and $k \geq p_0$ (the correctly specified case). However, the equivalence introduced by (4.9) seems to be more relevant in situations where the predictive capabilities of the two criteria after order selection are emphasized.

**5. Optimal prediction for possibly degenerate AR($\infty$) processes.** This section deals with optimal prediction problems in situations where the underlying AR($\infty$) process can degenerate to an AR process of finite order. We first adopt (K.5$'$) to replace the truly infinite-order assumption, (K.5).

(K.5$'$): The AR coefficients satisfy either

(i) $a_{p_0} \neq 0$ for some unknown $1 \leq p_0 < \infty$ and $a_l = 0$ for all $l \geq p_0 + 1$ or
(ii) (3.10).

From a practical point of view, (K.5$'$) is reasonably flexible because it contains any causal and invertible $\text{ARMA}(p,q)$ model, with $p + q \geq 1$, as a special case. Before tackling order selection problems under (K.5$'$), a preliminary result is needed, which shows that $\text{APE}_{\delta_n}$ and $\text{IC}_{P_n}$, with $\delta_n$ and $P_n$ satisfying certain conditions, are asymptotically efficient in finite-order cases.

THEOREM 4. *Assume that* (K.1)–(K.4) *and* (i) *of (*K.5$'$*) hold. Then* (2.8) *and* (2.9) *hold for* $(\hat{k}_{n,OS}, D_n) = (\hat{k}_{n,\delta_n}, D_{\text{APE}_{\delta_n}})$ *and* $(\hat{k}_{n,P_n}, P_n)$, *where* $\delta_n$ *satisfies* $\delta_n^{-1} \to \infty$ *and* (3.5), *and* $P_n$ *satisfies* $P_n \to \infty$ *and* $P_n = O(n^s)$ *for some* $0 < s < 1$. *In addition,* (2.3) [(2.4)] *is satisfied by these criteria.*

REMARK 9. Since Theorem 4 adopts $\{\text{AR}(1), \ldots, \text{AR}(K_n)\}$ as the set of candidate models, where $K_n \to \infty$ at a certain rate, the true model $\text{AR}(p_0)$ is included asymptotically. Zheng and Loh [31] also took this approach. However, unlike Theorem 4, their main concern is with the consistency in order selection.

When (ii) of (K.5′) holds, Example 6 points out that $\text{IC}_{P_n}$, with $P_n = o(\log n)$, possesses asymptotic efficiency. On the other hand, if (i) of (K.5′) is true, then Theorem 4 shows that $\text{IC}_{P_n}$, with $P_n \to \infty$ and $P_n = O(n^s)$, $0 < s < 1$, is asymptotically efficient under (K.1)–(K.4). These results taken together suggest that $\text{IC}_{P_n}$, with $P_n \to \infty$ and $P_n = o(\log n)$, simultaneously achieve (2.3) over the two types of AR processes defined in (i) and (ii) of (K.5′). According to Example 1 and Theorem 4, $\text{APE}_{\delta_n}$, with $\delta_n^{-1} \to \infty$ and $\log \delta_n^{-1} = o(\log n)$, also has this property. This discussion is now summarized in the following theorem.

THEOREM 5. *Assume that* (K.1)–(K.4) *and* (K.5′) *hold. Then* (2.3) [(2.4)] *holds for* $\hat{k}_n = \hat{k}_{n,\delta_n}$ *and* $\hat{k}_{n,P_n}$, *where* $\delta_n$ *satisfies* $\delta_n^{-1} \to \infty$ *and* $\log \delta_n^{-1} = o(\log n)$, *and* $P_n$ *satisfies* $P_n \to \infty$ *and* $P_n = o(\log n)$.

As pointed out in Examples 3 and 8, the criteria given by Theorem 5 fail to preserve asymptotic efficiency when (2.18) is included in (K.5′). To overcome this difficulty, we propose using an alternative criterion that chooses order $\hat{k}_n^{(\iota)}$:

$$(5.1) \qquad \hat{k}_n^{(\iota)} = \hat{k}_{n,2} I_{\{\hat{k}_{n,P_n} \neq \hat{k}_{n^\iota,P_{n^\iota}}\}} + \hat{k}_{n,P_n} I_{\{\hat{k}_{n,P_n} = \hat{k}_{n^\iota,P_{n^\iota}}\}},$$

where $0 < \iota < 1$, $P_n \to \infty$, $\hat{k}_{n^\iota,P_{n^\iota}} = \arg\min_{1 \leq k \leq K_{n^\iota}} \text{IC}_{P_{n^\iota}}(k)$ and

$$\text{IC}_{P_{n^\iota}}(k) = \log \hat{\sigma}_{n^\iota}^2(k) + \frac{P_{n^\iota} k}{n^\iota},$$

with $\hat{\sigma}_{n^\iota}^2(k) = (1/N_\iota) \sum_{j=K_{n^\iota}}^{n^\iota - 1} (x_{j+1} + \hat{\mathbf{a}}'_{n^\iota}(k)\mathbf{x}_j(k))^2$, $N_\iota = n^\iota - K_{n^\iota}$,

$$\hat{\mathbf{a}}_{n^\iota}(k) = -\hat{R}_{\iota,n^\iota}^{-1}(k)(1/N_\iota) \sum_{j=K_{n^\iota}}^{n^\iota - 1} \mathbf{x}_j(k) x_{j+1},$$

and $\hat{R}_{\iota,n^\iota}(k) = (1/N_\iota) \sum_{j=K_{n^\iota}}^{n^\iota - 1} \mathbf{x}_j(k)\mathbf{x}'_j(k)$ (note that without loss of generality, $n^\iota$ and $K_{n^\iota}$ are assumed to be positive integers). As observed, (5.1) is a hybrid selection procedure that combines together AIC and a BIC-like criterion. If the true order is finite, then it is expected that the orders selected by the BIC-like criterion at stages $n^\iota$ and $n$ will ultimately be the same due to consistency. On the other hand, when the true order is infinite, an interesting result is derived for which it is nearly impossible for the BIC-like criterion to choose the same order at these different stages; see Appendix D. Therefore, it is reasonable to adopt $\hat{k}_{n,2}$ (the order selected by AIC) if $\text{IC}_{P_n}$ and $\text{IC}_{P_{n^\iota}}$ determine different orders, and $\hat{k}_{n,P_n}$ (the order selected by the BIC-like criterion) otherwise. Theorem 6 justifies the validity of $\hat{k}_n^{(\iota)}$.

THEOREM 6. *Let* (K.1)–(K.4) *and* (K.6) *(with $D_n = 2$) hold, and $\iota$ and $P_n$ in (5.1) satisfy $0 < \iota < 1$, $P_n \to \infty$, $P_n = O(n^{\iota_1})$, with $0 < \iota_1 < (1 + \delta_1^*)/(2 + \delta_1^*)$, and $P_n/P_{n^\iota}^\nu = O(1)$ for some $\nu > 0$. Further, assume that the AR coefficients meet either of the following conditions:*

(i) (i) *of* (K.5′);

(ii) *for any $\xi > 0$, there are a nonnegative exponent, $0 \leq \theta = \theta(\xi) < 1$, and a positive number, $M = M(\xi) > 0$, such that*

$$(5.2) \qquad \liminf_{n \to \infty} \min_{k \in A_{P_n,\theta,M}} (k_{n,P_n}^*)^\xi \frac{L_{n,P_n}(k) - L_{n,P_n}(k_{n,P_n}^*)}{L_{n,P_n}(k_{n,P_n}^*)} > 0,$$

*and for all sufficiently large $n$,*

$$(5.3) \qquad A_{P_{n^\iota},\theta,M}^C \cap A_{P_n,\theta,M}^C = \varnothing,$$

*where $A_{P_n,\theta,M}$ is defined in* (K.6), $\varnothing$ *denotes the empty set,*

$$A_{P_n,\theta,M}^C = \{k : 1 \leq k \leq K_n, k \notin A_{P_n,\theta,M}\}$$

*and*

$$A_{P_{n^\iota},\theta,M}^C = \{k : 1 \leq k \leq K_{n^\iota}, k \notin A_{P_{n^\iota},\theta,M}\}.$$

*[Note that (5.3) implicitly implies that $a_l \neq 0$ for infinitely many $l$.]*

*Then* (2.3) *[*(2.4)*] holds for $\hat{k}_n = \hat{k}_n^{(\iota)}$.*

As an application of Theorem 6, it is shown in Example 9 below that $\hat{k}_n^{(\iota)}$, $0 < \iota < 1$, is asymptotically efficient when the true model is either (i) an AR process of finite order, (ii) an AR($\infty$) process with coefficients satisfying (3.10) (the exponential-decay case) or (iii) an AR($\infty$) process with coefficients satisfying (2.18) (the algebraic-decay case). To simplify the discussion, let $P_n$ be given by (4.6) with $C_1, C_2 > 0$, which satisfies all requirements for $P_n$ imposed by Theorem 6.

EXAMPLE 9. Assume that (K.1)–(K.4) hold, and either (K.5′) or (2.18) is satisfied. To show that $\hat{k}_n^{(\iota)}$, $0 < \iota < 1$, is asymptotically efficient in this situation, in view of Theorem 6 and the discussion after (K.6), it suffices to show that (5.2) and (5.3) are guaranteed by (3.10) as well as (2.18). First, assume that (3.10) is true. Then Example 6 shows that (3.12), with $D_n = P_n$, is valid, which further implies that for any $\xi > 0$, (5.2) holds for any $1 - \min\{\xi, 1\} < \theta < 1$ and any $M > 0$. In addition, (5.3) follows from (3.12) (with $D_n = P_n$), (4.6) (with $C_1, C_2 > 0$), $0 < \iota < 1$, and $0 < \theta < 1$. Next, let (2.18) hold. Reasoning as for Example 8, we obtain (3.20) (with $D_n = P_n$) and that for any $\xi > 0$, (5.2) holds for any $1 - (\min\{\xi, 2\}/2) < \theta < 1$ and any $M > 0$. Moreover, (5.3) is guaranteed by (3.20) (with $D_n = P_n$), (4.6) (with $C_1, C_2 > 0$), $0 < \iota < 1$, and $0 < \theta < 1$. Consequently, the desired result follows.

**6. Simulation results.** To illustrate the practical implications of our theoretical results, a simulation study is conducted in this section. Let observations be generated from the ARMA(1, 1) model

$$x_{t+1} = \phi_0 x_t + \varepsilon_{t+1} + \theta_0 \varepsilon_t,$$

where the $\varepsilon_t$'s are independent and identically $\mathcal{N}(0, 1)$ distributed and $(\phi_0, \theta_0) = (0.0, 0.98), (0.5, 0.8), (0.5, 0.4)$, and $(0.9, 0)$. For each combination of $(\phi_0, \theta_0)$, the empirical estimates of

$$RE(\hat{k}_n) = \frac{E(x_{n+1} - \hat{x}_{n+1}(\hat{k}_{n,2}))^2 - 1}{E(x_{n+1} - \hat{x}_{n+1}(\hat{k}_n))^2 - 1}$$

denoted by $\hat{RE}(\hat{k}_n)$ are obtained based on 5000 replications for $n = 180$, 300, 500, and 1000, where $\hat{k}_n = \hat{k}_{n,\delta_n}, \hat{k}_{n,P_n}, \hat{k}_n^{(\iota)}$, with $\delta_n = (\log n)^{-1}, 1 - (2/3)(\log n)^{-0.1}, 1 - (2/3)(\log n)^{-0.12}, 1 - (2/3)(\log n)^{-0.14}, P_n = 2.001 \log_2 n$ (HQ), $\log n$ (BIC), and $\iota = 0.69, 0.72, 0.75$. The penalty term of the BIC-like criterion associated with $\hat{k}_n^{(\iota)}$ is given by (4.6) with $C_1 = 0.8$ and $C_2 = 1$. In addition, $K_n$ and $K_{n^\iota}$ are set to the largest integers less than or equal to $n^{1/2}$ and $n^{\iota/2}$, respectively. Obviously, $RE(\hat{k}_n)$ measures the relative prediction efficiency of $\hat{x}_{n+1}(\hat{k}_n)$ to $\hat{x}_{n+1}(\hat{k}_{n,2})$ (AIC), and $RE(\hat{k}_n) > 1$ [$RE(\hat{k}_n) < 1$] suggests that $\hat{k}_n$ performs better (worse) than $\hat{k}_{n,2}$. These empirical results (see Table 1) are summarized as follows:

(1) AIC and BIC. The relative efficiencies of AIC and BIC are clearly affected by the magnitude of the MA parameter in finite-sample situations. Table 1 shows that when $\theta_0 \geq 0.8$, AIC notably outperforms BIC, which coincides with our theoretical findings in Examples 7 and 8 that BIC is not asymptotically efficient in truly AR($\infty$) models. In contrast, values of $\hat{RE}(\hat{k}_{n,\log n})$ are larger than 1 when $\theta_0 = 0.4$. However, since these values rapidly decrease from 1.26 to 1.08 as $n$ grows from 180 to 1000, the theoretical result just mentioned does not seem to be seriously violated. On the other hand, when $\theta_0 = 0$, values of $\hat{RE}(\hat{k}_{n,\log n})$ are larger than 3.5, and do not exhibit any decreasing trend. This matches the fact that BIC is consistent and asymptotically efficient in finite-order AR models (see Section 5), whereas AIC is not.

TABLE 1
*Empirical estimates of $RE(\hat{k}_n)$*

| $n$ | Models $(\phi_0, \theta_0)$ | $APE_{\delta_n}$ | | | | IC | | Two-stage | | |
|---|---|---|---|---|---|---|---|---|---|---|
| | | $\delta_{1,n}$ | $\delta_{2,n}$ | $\delta_{3,n}$ | $\delta_{4,n}$ | HQ | BIC | $n^{0.69}$ | $n^{0.72}$ | $n^{0.75}$ |
| 180 | (0.0, 0.98) | 0.88 | 0.93 | 0.92 | 0.93 | 0.89 | 0.78 | 0.95 | 0.94 | 0.94 |
| | (0.5, 0.8) | 0.95 | 0.95 | 0.95 | 0.94 | 0.98 | 0.83 | 0.98 | 0.97 | 0.97 |
| | (0.5, 0.4) | 1.28 | 1.07 | 1.05 | 1.03 | 1.36 | 1.26 | 1.08 | 1.08 | 1.08 |
| | (0.9, 0.0) | 2.21 | 1.34 | 1.33 | 1.28 | 2.31 | 3.59 | 1.81 | 1.86 | 1.95 |
| 300 | (0.0, 0.98) | 0.88 | 0.94 | 0.94 | 0.94 | 0.89 | 0.74 | 0.97 | 0.96 | 0.95 |
| | (0.5, 0.8) | 0.98 | 0.99 | 0.98 | 0.97 | 0.96 | 0.79 | 0.95 | 0.94 | 0.93 |
| | (0.5, 0.4) | 1.28 | 1.03 | 1.03 | 1.03 | 1.24 | 1.24 | 1.09 | 1.09 | 1.09 |
| | (0.9, 0.0) | 2.18 | 1.37 | 1.32 | 1.26 | 2.44 | 3.46 | 1.95 | 1.99 | 2.07 |
| 500 | (0.0, 0.98) | 0.85 | 0.94 | 0.95 | 0.95 | 0.85 | 0.68 | 0.96 | 0.95 | 0.94 |
| | (0.5, 0.8) | 0.97 | 0.97 | 0.97 | 0.96 | 0.98 | 0.78 | 0.97 | 0.95 | 0.95 |
| | (0.5, 0.4) | 1.28 | 1.10 | 1.05 | 1.04 | 1.32 | 1.17 | 1.03 | 1.02 | 1.06 |
| | (0.9, 0.0) | 2.31 | 1.36 | 1.31 | 1.27 | 2.64 | 4.17 | 2.39 | 2.43 | 2.41 |
| 1000 | (0.0, 0.98) | 0.86 | 0.95 | 0.96 | 0.95 | 0.86 | 0.66 | 0.99 | 0.98 | 0.98 |
| | (0.5, 0.8) | 1.05 | 0.97 | 0.96 | 0.96 | 1.01 | 0.80 | 0.97 | 0.97 | 0.95 |
| | (0.5, 0.4) | 1.36 | 1.12 | 1.09 | 1.04 | 1.37 | 1.08 | 1.00 | 1.00 | 0.98 |
| | (0.9, 0.0) | 2.33 | 1.27 | 1.26 | 1.21 | 2.86 | 4.07 | 2.65 | 2.74 | 2.67 |

Note: $\delta_{1,n} = (\log n)^{-1}, \delta_{2,n} = 1 - (2/3)(\log n)^{-0.1}, \delta_{3,n} = 1 - (2/3)(\log n)^{-0.12}$, and $\delta_{4,n} = 1 - (2/3)(\log n)^{-0.14}$.

(2) HQ and $APE_{\delta_{1,n}}$, where $\delta_{1,n} = (\log n)^{-1}$. First note that the prediction efficiencies of these two criteria seem quite close. They perform comparably to AIC when $\theta_0 = 0.8$, and much better than it when $\theta_0 \leq 0.4$. This phenomenon can be explained by the fact that HQ and $APE_{\delta_{1,n}}$ are asymptotically efficient in both the finite-order AR model and the $AR(\infty)$ model with AR coefficients decaying exponentially (see Theorem 5). Their efficiencies, however, are smaller than AIC in the case $\theta_0 = 0.98$. Since it is difficult to distinguish between an MA(1) process with a very large MA coefficient and an $AR(\infty)$ process with AR coefficients decaying algebraically in finite samples, Examples 3 and 8 (which show that HQ and $APE_{\delta_{1,n}}$ are not asymptotically efficient in the algebraic-decay case) may explain why HQ and $APE_{\delta_{1,n}}$ perform worse than AIC when $\theta_0$ is very close to unity. In addition, we also observe that these two criteria are not as efficient as BIC in the case $\theta_0 = 0$, but they beat BIC in all other cases.

(3) $APE_{\delta_{i,n}}, i = 2, 3, 4$, where $\delta_{2,n} = 1 - (2/3)(\log n)^{-0.1}, \delta_{3,n} = 1 - (2/3) \times (\log n)^{-0.12}$ and $\delta_{4,n} = 1 - (2/3)(\log n)^{-0.14}$. Table 1 shows that $APE_{\delta_{i,n}}, i = 2, 3, 4$, holds a slight advantage (disadvantage) over AIC when $\theta_0 = 0.4$ ($\theta_0 \geq 0.8$). However, since the amount of the advantage (disadvantage) is not sizable, these Monte Carlo results seem to support the theoretical findings revealed in Examples 4–6 and 8 that AIC and these $APE_{\delta_n}$ criteria

are asymptotically equivalent in both the exponential- and algebraic-decay cases. On the other hand, these criteria tend to outperform AIC when $\theta_0 = 0$. But they are still much less efficient than all other criteria in this case due to the lack of consistency in the finite-order case (see Remark 5).

(4) Two-stage criteria. One special feature of two-stage criteria is that they behave like AIC in situations where AIC dominates, and improve substantially over AIC in situations where AIC performs poorly. More specifically, values of $\hat{RE}(\hat{k}_n^{(\iota)})$ are rather close to 1 when $\theta_0 \geq 0.8$, and significantly larger than 1 when $\theta_0 = 0$. In the case $\theta_0 = 0.4$, the two-stage criteria perform slightly better than AIC when $n \leq 300$, and comparably to AIC when $n \geq 500$. These simulation results seem to match quite well with the conclusion drawn from Example 9 that the two-stage criteria are asymptotically efficient in all three (finite-order, exponential-decay, and algebraic-decay) cases. When $\theta_0 = 0$, the prediction performance of the two-stage criteria is similar to that of HQ and $\mathrm{APE}_{\delta_{1,n}}$ (particularly when $n$ is large), but worse than that of BIC.

In conclusion, note that the finite-sample behavior of the criteria considered in this section can be well predicted by the asymptotic results obtained in Sections 3–5. Some desirable features (when compared to AIC or BIC) of the $\mathrm{APE}_{\delta_{i,n}}$, HQ, and two-stage criteria are particularly encouraging (see the discussions above). The tuning parameters adopted in this section may also serve as good initial values for pursuing better finite-sample efficiencies.

## APPENDIX A: PROOFS OF PROPOSITION 2 AND (2.13)

PROOF OF PROPOSITION 2. In view of (1.7) and (2.1),

$$
\begin{aligned}
\frac{q_n(\hat{k}_{n,\mathrm{OS}}) - \sigma^2}{L_n(k_{n,D_n}^*)} \\
(\mathrm{A.1}) \quad &= E\{\mathbf{f}(\hat{k}_{n,\mathrm{OS}}) - \mathbf{f}(k_{n,D_n}^*) + \mathcal{S}(\hat{k}_{n,\mathrm{OS}}) \\
&\quad - \mathcal{S}(k_{n,D_n}^*) + \mathbf{f}(k_{n,D_n}^*) + \mathcal{S}(k_{n,D_n}^*)\}^2 (L_n(k_{n,D_n}^*))^{-1}.
\end{aligned}
$$

It is also not difficult to see that

$$
(\mathrm{A.2}) \qquad \frac{L_{n,D_n}(k_{n,D_n}^*)}{L_n(k_{n,D_n}^*)} = O(D_n - 1).
$$

By (A.1), (A.2), (2.2), (2.8) and (2.9), (2.10) follows. Moreover, if (2.11) is assumed, (2.10) can be rewritten as

$$
(\mathrm{A.3}) \qquad \lim_{n \to \infty} \frac{q_n(\hat{k}_{n,\mathrm{OS}}) - \sigma^2}{L_n(k_n^*)} = 1,
$$

and hence (2.3) [or (2.4)] holds for $\hat{k}_n = \hat{k}_{n,\text{OS}}$. □

PROOF OF (2.13). Under (K.1), (K.5), $\sup_{-\infty<t<\infty} E|e_t|^4 < \infty$ and $K_n = o(n^{1/2})$, an argument given in [14], page 2448 yields

(A.4) $$\operatorname*{p-lim}_{n\to\infty} \frac{E\{(y_{n+1} - \hat{y}_{n+1}(\hat{k}_{n,\text{OS}}))^2 | x_1, \ldots, x_n\} - \sigma^2}{L_n(\hat{k}_{n,\text{OS}})} = 1.$$

The desired result now follows from (2.12) and (A.4). □

## APPENDIX B: PROOFS OF THEOREMS 1 AND 2

In the rest of this paper, $C$ is used to denote a generic positive constant independent of the sample size $n$ and of any index with an upper (or lower) limit dependent on $n$. It also may have different values in different places. We start with a modification of Lemma 6 of [14].

LEMMA B.1. *Assume that* (K.1) *holds with* $\sum_{i=1}^{\infty} |i^{1/2} a_i| < \infty$ *replaced by* $\sum_{i=1}^{\infty} |a_i| < \infty$ *and* $\sup_{-\infty<t<\infty} E|e_t|^{2q} < \infty$ *for some* $q \geq 2$. *Let* $\{m_{i,n}\}$, $i = 0, 1, 2$, *be sequences of positive integers satisfying* $m_{2,n} \geq m_{1,n} \geq m_{0,n}$ *for all* $n \geq 1$. *Then, for all* $n \geq 1$ *and all* $1 \leq k, j \leq m_{0,n}$,

(B.1) $$E|S^2_{m_{1,n},m_{2,n}}(k) - \sigma_k^2 - (S^2_{m_{1,n},m_{2,n}}(j) - \sigma_j^2)|^q \leq C m_n^{-q/2} \|\mathbf{a}(j) - \mathbf{a}(k)\|_R^q,$$

*where* $m_n = m_{2,n} - m_{1,n} + 1$, $S^2_{m_{1,n},m_{2,n}}(k) = (1/m_n) \sum_{t=m_{1,n}}^{m_{2,n}} e^2_{t+1,k}$, $\mathbf{a}(j)$ *and* $\mathbf{a}(k)$ *in* (B.1) *are viewed as infinite-dimensional vectors with undefined entries set to zero and* $\sigma_k^2 = E(e^2_{1,k})$. *Also note that* $\|\mathbf{a}(j) - \mathbf{a}(k)\|_R^2 = |\|\mathbf{a} - \mathbf{a}(j)\|_R^2 - \|\mathbf{a} - \mathbf{a}(k)\|_R^2|$.

The proof of (B.1) is similar to that of [14], Lemma 6, and hence is omitted. Let $n\delta_n$, with $1/n \leq \delta_n \leq 1 - (1/n)$, be a positive integer. According to (3.3), for $k \neq k^*_{n,D_n}$,

(B.2) $$P(\hat{k}_{n,\delta_n} = k) \leq P\left(\frac{\text{APE}_{\delta_n}(k)}{U_{\delta_n}(k)} \leq \frac{\text{APE}_{\delta_n}(k^*_{n,D_n})}{U_{\delta_n}(k)}\right)$$
$$\leq \sum_{l=1}^{12} P\left(|N_l(k)| \geq \frac{1}{12} V_{n,D_n}(k)\right),$$

where $D_n = D_{\text{APE}_{\delta_n}}$ is used throughout this appendix, $U_{\delta_n}(k) = n(1-\delta_n)L_{n,D_n}(k)$,

$$V_{n,D_n}(k) = (L_{n,D_n}(k) - L_{n,D_n}(k^*_{n,D_n}))L^{-1}_{n,D_n}(k),$$

$$|N_1(k)| = U^{-1}_{\delta_n}(k)\left|\left(\sum_{i=n\delta_n}^{n-1} h_i(k)e^2_{i+1,k}\right) - k\sigma^2 \log \delta_n^{-1}\right|,$$

$$|N_2(k)| = U_{\delta_n}^{-1}(k) \left| \left( \sum_{i=n\delta_n}^{n-1} h_i(k_{n,D_n}^*) e_{i+1,k_{n,D_n}^*}^2 \right) - k_{n,D_n}^* \sigma^2 \log \delta_n^{-1} \right|,$$

$$|N_3(k)| = U_{\delta_n}^{-1}(k) |Q_{n\delta_n}(k) - k\sigma^2|,$$

$$|N_4(k)| = U_{\delta_n}^{-1}(k) |Q_{n\delta_n}(k_{n,D_n}^*) - k_{n,D_n}^* \sigma^2|,$$

$$|N_5(k)| = U_{\delta_n}^{-1}(k) |Q_n(k) - k\sigma^2|,$$

$$|N_6(k)| = U_{\delta_n}^{-1}(k) |Q_n(k_{n,D_n}^*) - k_{n,D_n}^* \sigma^2|,$$

$$|N_7(k)| = U_{\delta_n}^{-1}(k) \left| \sum_{i=n\delta_n}^{n-1} h_i(k) \hat{e}_{i,k}^2 \right|,$$

$$|N_8(k)| = U_{\delta_n}^{-1}(k) \left| \sum_{i=n\delta_n}^{n-1} h_i(k_{n,D_n}^*) \hat{e}_{i,k_{n,D_n}^*}^2 \right|,$$

$$|N_9(k)| = 2 U_{\delta_n}^{-1}(k) \left| \sum_{i=n\delta_n}^{n-1} h_i(k) \hat{e}_{i,k} e_{i+1,k} \right|,$$

$$|N_{10}(k)| = 2 U_{\delta_n}^{-1}(k) \left| \sum_{i=n\delta_n}^{n-1} h_i(k_{n,D_n}^*) \hat{e}_{i,k_{n,D_n}^*} e_{i+1,k_{n,D_n}^*} \right|,$$

$$|N_{11}(k)| = U_{\delta_n}^{-1}(k) \left| \sum_{i=n\delta_n}^{n-1} \{ \varepsilon_{i+1}^2(k) - \varepsilon_{i+1}^2(k_{n,D_n}^*) - \|\mathbf{a} - \mathbf{a}(k)\|_R^2 + \|\mathbf{a} - \mathbf{a}(k_{n,D_n}^*)\|_R^2 \} \right|,$$

$$|N_{12}(k)| = U_{\delta_n}^{-1}(k) \left| \sum_{i=n\delta_n}^{n-1} (e_{i+1,k} - e_{i+1,k_{n,D_n}^*}) e_{i+1} \right|,$$

and $\varepsilon_{i+1}(k) = e_{i+1,k} - e_{i+1}$.

By (B.2), Chebyshev's inequality, and moment bounds for $|N_i|, i = 1, \ldots, 12$ (see Lemmas B.2–B.4 below), an upper bound for $P(\hat{k}_{n,\delta_n} = k)$ can be obtained. This upper bound plays an important role in verifying Theorems 1 and 2.

LEMMA B.2. *Let the assumptions of Proposition 1 hold and $1/n \leq \delta_n \leq 1 - (1/n)$ satisfy (3.5). Then, for $q > 0$, all $1 \leq k \leq K_n$ and all sufficiently large $n$,*

(B.3) $\quad E(|N_1(k)|^q) \leq C U_{\delta_n}^{-q}(k) \left\{ \dfrac{k^{2q}}{(n\delta_n)^{q/2}} + (\log \delta_n^{-1})^q k^q \|\mathbf{a} - \mathbf{a}(k)\|_R^{2q} \right\}$

*and*

$$E(|N_2(k)|^q) \leq C U_{\delta_n}^{-q}(k) \left\{ \frac{k_{n,D_n}^{*2q}}{(n\delta_n)^{q/2}} + (\log \delta_n^{-1})^q k_{n,D_n}^{*q} \|\mathbf{a} - \mathbf{a}(k_{n,D_n}^*)\|_R^{2q} \right\}.$$
(B.4)

PROOF. We only prove (B.3) because the proof of (B.4) is similar. Define $D(i,k) = \mathbf{x}_i'(k) R^{-1}(k) \mathbf{x}_i(k)(i+1-K_n)^{-1}$ and $E(i,k) = k(i+1-K_n)^{-1}$. Then

$$U_{\delta_n}(k)|N_1(k)| \leq \left| \sum_{i=n\delta_n}^{n-1} (h_i(k) - D(i,k)) e_{i+1,k}^2 \right|$$

$$+ \left| \sum_{i=n\delta_n}^{n-1} D(i,k)(e_{i+1,k}^2 - e_{i+1}^2) \right| + \left| \sum_{i=n\delta_n}^{n-1} D(i,k)(e_{i+1}^2 - \sigma^2) \right|$$

(B.5)
$$+ \left| \sum_{i=n\delta_n}^{n-1} (D(i,k) - E(i,k)) \sigma^2 \right| + \sigma^2 \left| \sum_{i=n\delta_n}^{n-1} E(i,k) - k \log \delta_n^{-1} \right|$$

$$\equiv \mathrm{I}(k) + \mathrm{II}(k) + \mathrm{III}(k) + \mathrm{IV}(k) + \mathrm{V}(k).$$

By (3.5) and [14], Proposition 1, we have, for any $q > 0$, all $n\delta_n \leq i \leq n-1$, all $1 \leq k \leq K_n$ and all sufficiently large $n$,

(B.6) $$E\|\hat{R}_{i+1}^{-1}(k) - R^{-1}(k)\|^q \leq C \frac{k^q}{(i+1-K_n)^{q/2}}.$$

Using [28], Lemma 2, and Jensen's inequality, it follows that for any $r > 0$, all $n\delta_n \leq i \leq n-1$ and all $1 \leq k \leq K_n$,

(B.7) $$E(\|\mathbf{x}_i(k)\|^r) \leq C k^{r/2}$$

and

(B.8) $$E|e_{i+1,k}|^r \leq C.$$

According to (B.6)–(B.8), Minkowski's inequality and Hölder's inequality we have, for $q \geq 1$, all $1 \leq k \leq K_n$ and all sufficiently large $n$,

(B.9) $$E(\mathrm{I}(k))^q \leq \left( \sum_{i=n\delta_n}^{n-1} \|(h_i(k) - D(i,k)) e_{i+1,k}^2\|_q \right)^q \leq C k^{2q} (n\delta_n)^{-q/2},$$

where for a random variable $z$ and positive number $s$, $\|z\|_s = E(|z|^s)^{1/s}$.

To deal with $\mathrm{II}(k)$, note that the first moment bound theorem of Findley and Wei [4] and Jensen's inequality yield for any $r > 0$, all $K_n \leq i \leq n-1$ and all $1 \leq k \leq K_n$,

(B.10) $$E(|\mathbf{x}_i'(k) R^{-1}(k) \mathbf{x}_i(k) - k|^r) \leq C k^{r/2}.$$

Reasoning as for (B.8), we have, for any $r > 0$, all $K_n \leq i \leq n-1$ and all $1 \leq k \leq K_n$,

(B.11) $$E(|\varepsilon_{i+1}(k)|^r) \leq C\|\mathbf{a} - \mathbf{a}(k)\|_R^r.$$

(B.10), (B.11), [28], Lemma 2, and an argument similar to that used for obtaining (B.9) together imply that for $q \geq 2$ and all $1 \leq k \leq K_n$,

(B.12)
$$E(\text{II}(k))^q \leq C \left\{ E\left| \sum_{i=n\delta_n}^{n-1} D(i,k)\varepsilon_{i+1}^2(k) \right|^q + E\left| \sum_{i=n\delta_n}^{n-1} D(i,k)\varepsilon_{i+1}(k)e_{i+1} \right|^q \right\}$$
$$\leq C\{(\log \delta_n^{-1})^q k^q \|\mathbf{a} - \mathbf{a}(k)\|_R^{2q} + (n\delta_n)^{-q/2} k^q \|\mathbf{a} - \mathbf{a}(k)\|_R^q\}.$$

Similarly,

(B.13) $$E(\text{III}(k))^q \leq CE\left( \sum_{i=n\delta_n}^{n-1} D^2(i,k) \right)^{q/2} \leq Ck^q (n\delta_n)^{-q/2}$$

holds for $q \geq 2$ and all $1 \leq k \leq K_n$.

To deal with $\text{IV}(k)$, it can be shown by some algebraic manipulations that

$$\text{IV}(k) = \sigma^2 \left| \frac{T_{n-1}(k)}{N} - \frac{T_{n\delta_n - 1}(k)}{n\delta_n - K_n} + \sum_{i=n\delta_n}^{n-1} \frac{T_{i-1}(k)}{(i - K_n)(i + 1 - K_n)} \right|,$$

where $T_i(k) = \sum_{j=K_n}^{i} \mathbf{x}_i'(k) R^{-1}(k) \mathbf{x}_i(k) - k$. By an argument similar to that given in the proof of Lemma 3 of [14] and Jensen's inequality, one has for any $q > 0$, all $n\delta_n - 1 \leq i \leq n-1$ and all $1 \leq k \leq K_n$,

$$E\left| \frac{T_i(k)}{i+1-K_n} \right|^q \leq C \frac{k^{3q/2}}{(i+1-K_n)^{q/2}}.$$

This and the Minkowski inequality yield that for $q \geq 1$ and all $1 \leq k \leq K_n$,

(B.14) $$E(\text{IV}(k))^q \leq Ck^{3q/2}(n\delta_n)^{-q/2}.$$

In addition, it is straightforward to show that for all $1 \leq k \leq K_n$,

(B.15) $$\text{V}(k) \leq C\left(\frac{1-\delta_n}{\delta_n}\right) k\left(\frac{K_n}{n}\right).$$

Consequently, (B.3) follows from (B.5), (B.9), (B.12)–(B.15), Jensen's inequality, and the fact that for any $r > 0$,

(B.16) $$\lim_{k \to \infty} k^r \|\mathbf{a} - \mathbf{a}(k)\|_R^{2r} = 0,$$

which is ensured by (K.1). □

LEMMA B.3. *Under the assumptions of Lemma B.2, for $q > 0$, all $1 \leq k \leq K_n$ and all sufficiently large $n$,*

(B.17) $\qquad E(|N_3(k)|^q) \leq C U_{\delta_n}^{-q}(k)\{k^{2q}(n\delta_n)^{-q/2} + k^{q/2}\},$

(B.18) $\qquad E(|N_4(k)|^q) \leq C U_{\delta_n}^{-q}(k)\{k_{n,D_n}^{*2q}(n\delta_n)^{-q/2} + k_{n,D_n}^{*q/2}\},$

(B.19) $\qquad E(|N_5(k)|^q) \leq C U_{\delta_n}^{-q}(k)(k^{2q}n^{-q/2} + k^{q/2}),$

(B.20) $\qquad E(|N_6(k)|^q) \leq C U_{\delta_n}^{-q}(k)(k_{n,D_n}^{*2q}n^{-q/2} + k_{n,D_n}^{*q/2}),$

(B.21) $\qquad E(|N_7(k)|^q) \leq C U_{\delta_n}^{-q}(k)k^{3q}(n\delta_n)^{-q},$

(B.22) $\qquad E(|N_8(k)|^q) \leq C U_{\delta_n}^{-q}(k)k_{n,D_n}^{*3q}(n\delta_n)^{-q},$

(B.23) $\qquad E(|N_9(k)|^q) \leq C U_{\delta_n}^{-q}(k)k^{2q}(n\delta_n)^{-q/2} \quad \text{and}$

(B.24) $\qquad E(|N_{10}(k)|^q) \leq C U_{\delta_n}^{-q}(k)k_{n,D_n}^{*2q}(n\delta_n)^{-q/2}.$

PROOF. See [12], Lemmas A.7 and A.8. □

LEMMA B.4. *Let the assumptions of Lemma B.1 hold and $1/n \leq \delta_n \leq 1 - (1/n)$. Then, for $q \geq 2$ and all $1 \leq k \leq K_n$, with $K_n \leq n\delta_n$,*

(B.25) $\qquad E(|N_i(k)|^q) \leq C U_{\delta_n}^{-q/2}(k) L_{n,D_n}^{-q/2}(k) \|\mathbf{a}(k) - \mathbf{a}(k_{n,\delta_n}^*)\|_R^q,$

*where $i = 11$ and $12$.*

PROOF. First note that

$$E|L_{n,D_n}(k)N_{11}(k)|^q$$

(B.26) $\qquad \leq E\left|\dfrac{\sum_{i=n\delta_n}^{n-1}\{e_{i+1,k}^2 - e_{i+1,k_{n,D_n}^*}^2 - \sigma_k^2 + \sigma_{k_{n,D_n}^*}^2\}}{n(1-\delta_n)}\right|^q$

$\qquad + E\left|\dfrac{\sum_{i=n\delta_n}^{n-1}(e_{i+1,k} - e_{i+1,k_{n,D_n}^*})e_{i+1}}{n(1-\delta_n)}\right|^q \equiv \text{(I)} + \text{(II)}.$

According to (B.1), one has for all $1 \leq k \leq K_n$,

(B.27) $\qquad \text{(I)} \leq \dfrac{C\|\mathbf{a}(k) - \mathbf{a}(k_{n,D_n}^*)\|_R^q}{(1-\delta_n)^{q/2}n^{q/2}}.$

Lemma 2 of [28] and the convexity of $x^{q/2}, x > 0$, yield for all $1 \leq k \leq K_n$,

(B.28) $\qquad \text{(II)} \leq \dfrac{C}{\{n(1-\delta_n)\}^{(q/2)+1}} \sum_{i=n\delta_n}^{n-1} E(|e_{i+1,k} - e_{i+1,k_{n,D_n}^*}|^q)$

$\qquad\qquad \leq \dfrac{C\|\mathbf{a}(k) - \mathbf{a}(k_{n,D_n}^*)\|_R^q}{(1-\delta_n)^{q/2}n^{q/2}}.$

Consequently, (B.25), with $i = 11$, is ensured by (B.26)–(B.28). The proof is completed by noting that (B.25), with $i = 12$, is an immediate consequence of (B.28). □

Armed with Lemmas B.2–B.4, we have the following result.

COROLLARY B.1. *Let* (K.1)–(K.5), (3.4) *and* (3.5) *hold. Then, for any* $r > 0$,

$$\lim_{n \to \infty} E\left(\frac{L_{n,D_n}(\hat{k}_{n,\delta_n})}{L_{n,D_n}(k^*_{n,D_n})} - 1\right)^r = 0. \tag{B.29}$$

PROOF. Define $I_{1,n}(k) = L_{n,D_n}(k)/L_{n,D_n}(k^*_{n,D_n})$. Let $\epsilon > 0$ be arbitrarily given. Then, by (B.2), one has

$$\begin{aligned}
&E(I_{1,n}(\hat{k}_{n,\delta_n}) - 1)^r \\
&= \sum_{k=1}^{K_n} (I_{1,n}(k) - 1)^r P(\hat{k}_{n,\delta_n} = k) \\
&\leq \epsilon^r + \sum_{l=1}^{12} \left\{ \sum_{k \in A^{(\delta_n)}_{\epsilon,n}} (I_{1,n}(k) - 1)^r P(|N_l(k)| \geq (1/12) V_{n,D_n}(k)) \right\},
\end{aligned} \tag{B.30}$$

where $A^{(\delta_n)}_{\epsilon,n} = \{k : 1 \leq k \leq K_n, I_{1,n}(k) - 1 > \epsilon\}$. In view of (B.30), (B.29) holds if for $l = 1, \ldots, 12$,

$$\lim_{n \to \infty} \sum_{k \in A^{(\delta_n)}_{\epsilon,n}} (I_{1,n}(k) - 1)^r P(|N_l(k)| > (1/12) V_{n,D_n}(k)) = 0. \tag{B.31}$$

In the following, we only prove (B.31) for $l = 1, 3$ and 11 because the proofs for $l = 2, 7, 8, 9$ and 10 are similar to that for $l = 1$, proofs for $l = 4, 5$ and 6 are similar to that for $l = 3$, and the proof for $l = 12$ is similar to that for $l = 11$.

By (B.3), Chebyshev's inequality, (3.4), (3.5) and the facts that

$$L_{n,D_n}(k) \geq \|\mathbf{a} - \mathbf{a}(k)\|_R^2, \qquad nL_{n,D_n}(k) \geq C\frac{k \log \delta_n^{-1}}{1 - \delta_n}, \tag{B.32}$$

and $I_{1,n}(k) \leq C/L_{n,D_n}(k^*_{n,D_n})$ if $1 \leq k \leq k^*_{n,D_n}$ and $I_{1,n}(k) \leq Ck/k^*_{n,D_n}$ if $k^*_{n,D_n} < k \leq K_n$, we have, for sufficiently large $q$,

$$\begin{aligned}
&\sum_{k \in A^{(\delta_n)}_{\epsilon,n}} (I_{1,n}(k) - 1)^r P(|N_1(k)| > (1/12) V_{n,D_n}(k)) \\
&\leq C \sum_{k \in A^{(\delta_n)}_{\epsilon,n}} I^r_{1,n}(k) V^{-(q-r)}_{n,D_n}(k) \{k^q f^{-q}_{1,n} + (f_{2,n}k)^q n^{-q}\}
\end{aligned}$$

(B.33)
$$\leq C\left(\frac{1+\epsilon}{\epsilon}\right)^{q-r}\left[\sum_{k=1}^{k^*_{n,D_n}}\left\{\frac{k^q(1-\delta_n)^r}{f^q_{1,n}(k^*_{n,D_n}f_{3,n}\delta_n)^r}+\frac{f^{q-r}_{2,n}k^q}{k^{*r}_{n,D_n}n^{q-r}}\right\}\right.$$

$$\left.+\sum_{k=k^*_{n,D_n}+1}^{K_n}\left\{\frac{k^{q+r}}{k^{*r}_{n,D_n}f^q_{1,n}}+\frac{f^q_{2,n}k^{q+r}}{k^{*r}_{n,D_n}n^q}\right\}\right]=o(1),$$

where $f_{1,n}=(\log\delta_n^{-1})(n\delta_n)^{1/2}$, $f_{2,n}=\log\delta_n^{-1}/(1-\delta_n)$ and $f_{3,n}=\log\delta_n^{-1}/(n\delta_n)$. Therefore, (B.31) holds for $l=1$.

By (B.17), (B.25), an argument similar to that used for obtaining (B.33) and the fact that $k^*_{n,D_n}\to\infty$ as $n\to\infty$, for sufficiently large $q$,

$$\sum_{k\in A^{(\delta_n)}_{\epsilon,n}}(I_{1,n}(k)-1)^r P(|N_3(k)|>(1/12)V_{n,D_n}(k))$$

(B.34)
$$\leq C\left(\frac{1+\epsilon}{\epsilon}\right)^{q-r}\left\{\sum_{k=1}^{K_n}I^r_{1,n}(k)k^q f^{-q}_{1,n}+\sum_{k=1}^{k^*_{n,D_n}}k^{q/2}U^{-q}_{\delta_n}(k^*_{n,D_n})\right.$$

$$\left.+\sum_{k=k^*_{n,D_n}+1}^{K_n}k^{q/2}U^{-q}_{\delta_n}(k)I^r_{1,n}(k)\right\}=o(1)$$

and

$$\sum_{k\in A^{(\delta_n)}_{\epsilon,n}}(I_{1,n}(k)-1)^r P(|N_{11}(k)|>(1/12)V_{n,D_n}(k))$$

(B.35)
$$\leq C\left(\frac{1+\epsilon}{\epsilon}\right)^{q-r}\left\{\sum_{k=1}^{k^*_{n,D_n}}U^{-q/2}_{\delta_n}(k^*_{n,D_n})+\sum_{k=k^*_{n,D_n}+1}^{K_n}U^{-q/2}_{\delta_n}(k)I^r_{1,n}(k)\right\}$$

$$=o(1).$$

In view of (B.33)–(B.35), the proof is complete. □

COROLLARY B.2. *Assume that* (K.1)–(K.6) *hold and $\delta_n$ satisfies* (3.4) *and* (3.5). *Then, for sufficiently large $q$,*

(B.36)
$$E\left|\frac{\mathcal{S}(\hat{k}_{n,\delta_n})-\mathcal{S}(k^*_{n,D_n})}{(L_{n,D_n}(\hat{k}_{n,\delta_n}))^{1/2}}\right|^{2q}=O((k^*_{n,D_n})^{-(1-\theta)q+\theta})+o((\log\delta_n^{-1})^{-q})$$

$$+O((\log\delta_n^{-1})^{-q/2}(k^*_{n,D_n})^{(-q/2)+\theta}),$$

*where $\mathcal{S}(k)$ is defined in Section* 2 *and $0\leq\theta=\theta(\xi)<1$ is any exponent obtained from* (K.6) *with $\xi$ satisfying* (3.7).

PROOF. Let $\xi$ satisfy (3.7). Then (K.6) guarantees that there are $0 \leq \theta = \theta(\xi) < 1$ and $M = M(\xi) > 0$ such that (2.16) is satisfied. Let $(\theta, M)$ be any such pair. Define $I_{2,n}(k) = L_{n,D_n}(k) - L_{n,D_n}(k^*_{n,D_n})$. By Hölder's inequality and the fact that for any $h > 0$,

$$
\begin{aligned}
E|\mathcal{S}(k) - \mathcal{S}(k^*_{n,D_n})|^{2h} &\leq C\|\mathbf{a}(k) - \mathbf{a}(k^*_{n,D_n})\|_R^{2h} \\
&\leq C(I_{2,n}(k) + f_{2,n}N^{-1}|k - k^*_{n,D_n}|\sigma^2)^h
\end{aligned}
\tag{B.37}
$$

(which follows from [28], Lemma 2, (K.3) and the definition of $k^*_{n,D_n}$), one has for $q > 0$ and $1 < r < \infty$,

$$
\begin{aligned}
E\left|\frac{\mathcal{S}(\hat{k}_{n,\delta_n}) - \mathcal{S}(k^*_{n,D_n})}{L^{1/2}_{n,D_n}(\hat{k}_{n,\delta_n})}\right|^{2q} \\
\leq \sum_{k=1}^{K_n} \left(E\left|\frac{\mathcal{S}(k) - \mathcal{S}(k^*_{n,D_n})}{L^{1/2}_{n,D_n}(k)}\right|^{2qr}\right)^{1/r} P^{(r-1)/r}(\hat{k}_{n,\delta_n} = k) \\
\leq C \sum_{k=1}^{K_n} \left\{V^q_{n,D_n}(k) + \left|\frac{f_{2,n}(k - k^*_{n,D_n})}{NL_{n,D_n}(k)}\right|^q\right\} P^{(r-1)/r}(\hat{k}_{n,\delta_n} = k) \\
\leq C\left\{\sum_{k=1}^{K_n} V^q_{n,D_n}(k) P^{(r-1)/r}(\hat{k}_{n,\delta_n} = k) + \sum_{\substack{k=1 \\ k \notin A_{D_n,\theta,M}}}^{K_n} \left|\frac{f_{2,n}(k - k^*_{n,D_n})}{NL_{n,D_n}(k)}\right|^q \right. \\
\left. + \sum_{\substack{k=1 \\ k \in A_{D_n,\theta,M}}}^{K_n} \left|\frac{f_{2,n}(k - k^*_{n,D_n})}{NL_{n,D_n}(k)}\right|^q P^{(r-1)/r}(\hat{k}_{n,\delta_n} = k)\right\} \\
\equiv C\{(\mathrm{I}) + (\mathrm{II}) + (\mathrm{III})\},
\end{aligned}
\tag{B.38}
$$

where $A_{D_n,\theta,M}$ is a set of positive integers defined in (K.6).

By the definitions of $A_{D_n,\theta,M}$, $L_{n,D_n}(k)$ and $L_{n,D_n}(k^*_{n,D_n})$, it is easy to see that

$$(\mathrm{II}) \leq C(k^*_{n,D_n})^{-(1-\theta)q+\theta}. \tag{B.39}$$

In view of (B.2) and the fact that for $a, b \geq 0, (a+b)^{(r-1)/r} \leq a^{(r-1)/r} + b^{(r-1)/r}$, one obtains

$$(\mathrm{I}) \leq \sum_{l=1}^{12}\left\{\sum_{k=1}^{K_n} V^q_{n,D_n}(k) P^{(r-1)/r}(|N_l(k)| \geq (1/12)V_{n,D_n}(k))\right\}. \tag{B.40}$$

In the following, we shall show that when $q$ is sufficiently large,

$$\sum_{k=1}^{K_n} V^q_{n,D_n}(k) P^{(r-1)/r}(|N_l(k)| \geq (1/12)V_{n,D_n}(k)) = o((\log \delta_n^{-1})^{-q}), \tag{B.41}$$

for $l = 1, \ldots, 10$; and

(B.42)
$$\sum_{k=1}^{K_n} V_{n,D_n}^q(k) P^{(r-1)/r}(|N_l(k)| \geq (1/12)V_{n,D_n}(k))$$
$$= O((\log \delta_n^{-1})^{-q/2}(k_{n,D_n}^*)^{(-q/2)+\theta}) + o((\log \delta_n^{-1})^{-q}),$$

for $l = 11$ and $12$. As a result, one has for sufficiently large $q$,

(B.43) $\quad (I) = O((\log \delta_n^{-1})^{-q/2}(k_{n,D_n}^*)^{(-q/2)+\theta}) + o((\log \delta_n^{-1})^{-q}).$

By Lemma B.2, (3.4), (3.5), (K.4) and (B.32), for sufficiently large $q$,

(B.44)
$$\sum_{k=1}^{K_n} V_{n,D_n}^q(k) P^{(r-1)/r}(|N_1(k)| \geq (1/12)V_{n,D_n}(k))$$
$$\leq C \sum_{k=1}^{K_n} [E\{|N_1(k)|^{qr/(r-1)}\}]^{(r-1)/r}$$
$$\leq \frac{C}{(\log \delta_n^{-1})^q} \left( \sum_{k=1}^{K_n} \frac{k^q}{(n\delta_n)^{q/2}} + \frac{(\log \delta_n^{-1})^{2q} k^q}{(1-\delta_n)^q n^q} \right) = o((\log \delta_n^{-1})^{-q}),$$

which yields (B.41) for $l = 1$. For $l = 3$, according to (3.5), Lemma B.3, (B.32) and the fact that $k_{n,D_n}^* \to \infty$ as $n \to \infty$, one has for sufficiently large $q$,

(B.45)
$$\sum_{k=1}^{K_n} V_{n,D_n}^q(k) P^{(r-1)/r}(|N_3(k)| \geq (1/12)V_{n,D_n}(k))$$
$$\leq C \sum_{k=1}^{K_n} [E\{|N_3(k)|^{qr/(r-1)}\}]^{(r-1)/r}$$
$$\leq \frac{C}{(\log \delta_n^{-1})^q} \left( \sum_{k=1}^{K_n} \frac{k^q}{(n\delta_n)^{q/2}} + \sum_{k=1}^{k_{n,D_n}^*} \frac{k^{q/2}}{k_{n,D_n}^{*q}} + \sum_{k=k_{n,D_n}^*+1}^{K_n} k^{-q/2} \right)$$
$$= o((\log \delta_n^{-1})^{-q}).$$

The proofs of (B.41) for $l = 2, 7, 8, 9$ and $10$ are similar to that of (B.44) and the proofs of (B.41) for $l = 4, 5$ and $6$ are similar to that of (B.45). We skip the details in order to save space. The proof of (B.42) is a bit more complicated. By (2.16), Lemma B.4, (B.37), (3.4) and the restriction on $\xi$, one has for sufficiently large $q$,

$$\sum_{k=1}^{K_n} V_{n,D_n}^q(k) P^{(r-1)/r}(|N_l(k)| \geq (1/12)V_{n,D_n}(k))$$

$$
\begin{aligned}
&\leq \sum_{\substack{k=1\\k\notin A_{D_n,\theta,M}}}^{K_n} \{E|N_{11}(k)|^{qr/(r-1)}\}^{(r-1)/r}\\
&\quad + \sum_{\substack{k=1\\k\in A_{D_n,\theta,M}}}^{K_n} V_{n,D_n}^{-q}(k)\{E|N_{11}(k)|^{2qr/(r-1)}\}^{(r-1)/r}
\end{aligned}
$$
(B.46)

$$
\begin{aligned}
&\leq C\Bigg\{ \sum_{\substack{k=1\\k\notin A_{D_n,\theta,M}}}^{K_n} \frac{\|\mathbf{a}(k)-\mathbf{a}(k_{n,D_n}^*)\|_R^q}{U_{\delta_n}^{q/2}(k)L_{n,D_n}^{q/2}(k)}\\
&\quad + \sum_{\substack{k=1\\k\in A_{D_n,\theta,M}}}^{K_n} \frac{I_{2,n}^q(k)+|(f_{2,n}(k-k_{n,D_n}^*))/N|^q}{U_{\delta_n}^q(k)I_{2,n}^q(k)}\Bigg\}\\
&\leq C\Bigg[\frac{k_{n,D_n}^{*\theta}}{(\log\delta_n^{-1}k_{n,D_n}^*)^{q/2}}\\
&\quad + \frac{\{1+(k_{n,D_n}^*)^{\xi q}\}}{(\log\delta_n^{-1})^q}\Bigg\{\sum_{k=1}^{k_{n,D_n}^*} k_{n,D_n}^{*-q} + \sum_{k=k_{n,D_n}^*+1}^{K_n} k^{-q}\Bigg\}\Bigg]\\
&= O((\log\delta_n^{-1})^{-q/2}(k_{n,D_n}^*)^{(-q/2)+\theta}) + o((\log\delta_n^{-1})^{-q}),
\end{aligned}
$$

where $l=11$ or $12$.

Following arguments similar to those used to obtain (B.40) and (B.44)–(B.46), it can be shown that

$$
\begin{aligned}
\text{(III)} &\leq \sum_{l=1}^{12}\Bigg\{\sum_{\substack{k=1\\k\in A_{D_n,\theta,M}}}^{K_n} \left|\frac{f_{2,n}(k-k_{n,D_n}^*)}{(NL_{n,D_n}(k))}\right|^q\\
&\quad \times P^{(r-1)/r}(|N_l(k)|\geq (1/12)V_{n,D_n}(k))\Bigg\}\\
&= o((\log\delta_n^{-1})^{-q}),
\end{aligned}
$$
(B.47)

where the equality holds for sufficiently large $q$. (For a detailed proof of (B.47), see [12], Corollary A.2.) Consequently, (B.36) is ensured by (B.38), (B.39), (B.43) and (B.47). $\square$

COROLLARY B.3. *Assume that the assumptions of Corollary B.2 hold. Then, for sufficiently large q,*

$$\text{(B.48)} \qquad \lim_{n\to\infty} E\left|\frac{\mathbf{f}(\hat{k}_{n,\delta_n}) - \mathbf{f}(k^*_{n,D_n})}{(L_{n,D_n}(\hat{k}_{n,\delta_n}))^{1/2}}\right|^{2q} = o((\log \delta_n^{-1})^q).$$

PROOF. Equation (B.48) can be verified using arguments similar to those in the proofs of [14], Lemmas 7 and 8, and Corollary B.2 above. For details, see [12], Corollary A.3. □

We are now ready to prove Theorem 1.

PROOF OF THEOREM 1. Let $(\eta, \theta)$ be a pair satisfying (3.6), where $\eta > 0$ and $0 \leq \theta = \theta(\xi) < 1$ is obtained from (K.6) with $\xi$ obeying (3.7). Then, by using Hölder's inequality, Jensen's inequality, Corollaries B.1 and B.2, (3.4) and (3.6), and taking $q > \max\{\eta^{-1}, 1\}$,

$$\text{(B.49)} \quad \begin{aligned} (D_n - 1)&E\left[\frac{\{\mathcal{S}(\hat{k}_{n,\delta_n}) - \mathcal{S}(k^*_{n,D_n})\}^2}{L_{n,D_n}(k^*_{n,D_n})}\right] \\ &\leq (D_n - 1)\left[E\frac{|\mathcal{S}(\hat{k}_{n,\delta_n}) - \mathcal{S}(k^*_{n,D_n})|^{2q}}{(L_{n,D_n}(\hat{k}_{n,\delta_n}))^q}\right]^{1/q} \\ &\quad \times \left[E\left\{\frac{L_{n,D_n}(\hat{k}_{n,\delta_n})}{L_{n,D_n}(k^*_{n,D_n})}\right\}^{q/(q-1)}\right]^{(q-1)/q} \\ &= O\left(\frac{D_n - 1}{(k^*_{n,D_n})^{1-\theta(1+1/q)}}\right) + o\left(\frac{D_n - 1}{\log \delta_n^{-1}}\right) \\ &\quad + O\left(\frac{D_n - 1}{(\log \delta_n^{-1})^{1/2}(k^*_{n,D_n})^{1/2-\theta/q}}\right) = o(1). \end{aligned}$$

By Corollaries B.1 and B.3 and an argument similar to that used to prove (B.49),

$$\text{(B.50)} \qquad (D_n - 1)E\left[\frac{\{\mathbf{f}(\hat{k}_{n,\delta_n}) - \mathbf{f}(k^*_{n,D_n})\}^2}{L_{n,D_n}(k^*_{n,D_n})}\right] = o(1).$$

Consequently, the desired result is ensured by (B.49), (B.50) and Proposition 2. □

PROOF OF THEOREM 2. First note that when $\lim_{n\to\infty} \delta_n = 1$ and condition (i) [or (ii)] of Theorem 2 are assumed instead of (3.4) and (3.5), the left-hand sides of (B.33)–(B.35) still converge to 0. Therefore, (B.29) follows. Let $0 < \xi < (1/2) - \xi_2$ if condition (i) of Theorem 2 holds, and $0 < \xi < \min\{(1/2) - \xi_2, (\delta_1^*/2) - \xi_2\}$ if condition (ii) of Theorem 2 holds. Then, by Jensen's inequality and the same reasoning used in the proofs of Corollaries B.2 and B.3, we have for any $q > 0$,

(B.51)
$$E\left|\frac{\mathcal{S}(\hat{k}_{n,\delta_n}) - \mathcal{S}(k_{n,D_n}^*)}{(L_{n,D_n}(\hat{k}_{n,\delta_n}))^{1/2}}\right|^{2q} = o(1) \quad \text{and}$$
$$E\left|\frac{\mathbf{f}(\hat{k}_{n,\delta_n}) - \mathbf{f}(k_{n,D_n}^*)}{(L_{n,D_n}(\hat{k}_{n,\delta_n}))^{1/2}}\right|^{2q} = o(1).$$

Consequently, the claimed result follows from (B.29), (B.51), (2.11) and Proposition 2. □

## APPENDIX C: PROOF OF THEOREM 3

Instead of verifying Theorem 3 directly, we will first investigate the prediction performance of $S_n^{(P_n)}(k)$, defined in (2.15). By analogy with (4.1) of [25],

(C.1)
$$\begin{aligned}S_n^{(P_n)}(k) &= NL_{n,D_n}(k) + P_nk(\hat{\sigma}_n^2(k) - \sigma^2) \\ &\quad + (k\sigma^2 - N\|\hat{\mathbf{a}}_n(k) - \mathbf{a}(k)\|_{\hat{R}_n(k)}^2) \\ &\quad + N\sigma^2 + N(S_{K_n,n-1}^2(k) - \sigma_k^2),\end{aligned}$$

where $D_n = P_n$, for a $k \times k$ symmetric matrix $A$ and a $k$-dimensional vector $\mathbf{y}$, $\|\mathbf{y}\|_A = \mathbf{y}'A\mathbf{y}$, and the definition of $S_{K_n,n-1}^2(k)$ can be found in Lemma B.1. Note that the relation, $D_n = P_n$, will be used throughout this appendix. Based on (C.1) and an argument similar to that used in (5.34) of [14], we have

(C.2) $$P(\hat{k}_{n,P_n}^S = k) \leq \sum_{i=1}^{5} P(|U_{in}(k)| \geq (1/5)V_{n,D_n}(k)),$$

where $V_{n,D_n}(k)$ is defined after (B.2), $\hat{k}_{n,P_n}^S = \arg\min_{1 \leq k \leq K_n} S_n^{(P_n)}(k)$,

$NL_{n,D_n}(k)|U_{1,n}(k)| = |P_nk(\hat{\sigma}_n^2(k) - \sigma^2)|,$

$NL_{n,D_n}(k)|U_{2,n}(k)| = |P_nk_{n,D_n}^*(\hat{\sigma}_n^2(k_{n,D_n}^*) - \sigma^2)|,$

$NL_{n,D_n}(k)|U_{3,n}(k)| = |k\sigma^2 - N\|\hat{\mathbf{a}}_n(k) - \mathbf{a}(k)\|_{\hat{R}_n(k)}^2|,$

$NL_{n,D_n}(k)|U_{4,n}(k)| = |k_{n,D_n}^*\sigma^2 - N\|\hat{\mathbf{a}}_n(k_{n,D_n}^*) - \mathbf{a}(k_{n,D_n}^*)\|_{\hat{R}_n(k_{n,D_n}^*)}^2|$ and

$NL_{n,D_n}(k)|U_{5,n}(k)| = |S_{K_n,n-1}^2(k) - \sigma_k^2 - S_{K_n,n-1}^2(k_{n,D_n}^*) - \sigma_{k_{n,D_n}^*}^2|.$

THEOREM C.1. *Let the assumptions of Theorem 3 hold. Then* (2.8) *and* (2.9) *hold for* $\hat{k}_{n,\mathrm{OS}} = \hat{k}_{n,P_n}^S$, *and* (4.4) *holds with* $\hat{k}_{n,P_n}$ *replaced by* $\hat{k}_{n,P_n}^S$.

PROOF. By Lemma B.1 and analogies with [14], (5.43) and (5.47), we have for $q > 0$, all $1 \leq k \leq K_n$ and all sufficiently large $n$, $E|U_{1,n}(k)|^q \leq C(P_n^q k^q N^{-q} + N^{-q/2})$, $E|U_{2,n}(k)|^q \leq C(P_n^q k_{n,D_n}^{*q} N^{-q} + N^{-q/2})$, $E|U_{3,n}(k)|^q \leq C(k^{2q}N^{-q/2} + k^{q/2})N^{-q}L_{n,D_n}^{-q}(k)$, $E|U_{4,n}(k)|^q \leq C(k_{n,D_n}^{*2q} N^{-q/2} + k_{n,D_n}^{*q/2})N^{-q} \times L_{n,D_n}^{-q}(k)$ and $E|U_{5,n}(k)|^q \leq C\|\mathbf{a}(k) - \mathbf{a}(k_{n,D_n}^*)\|_R^q N^{-q/2}L_{n,D_n}^{-q}(k)$. These moment bounds and an argument similar to that used to verify Corollary B.1 give for $q > 0$,

$$\lim_{n \to \infty} E\left(\frac{L_{n,D_n}(\hat{k}_{n,P_n}^S)}{L_{n,D_n}(k_{n,D_n}^*)} - 1\right)^q = 0. \tag{C.3}$$

With the help of (K.1)–(K.6), (4.1), (4.2) and the above moment properties, we can follow the ideas used in the proofs of Corollaries B.2 and B.3 to obtain that for $q > 0$,

$$E\left|\frac{\mathcal{S}(\hat{k}_{n,P_n}^S) - \mathcal{S}(k_{n,D_n}^*)}{(L_{n,D_n}(\hat{k}_{n,P_n}^S))^{1/2}}\right|^{2q} = O((k_{n,D_n}^*)^{-(1-\theta)q+\theta}) + o((P_n - 1)^{-q}) + O((P_n - 1)^{-q/2}(k_{n,D_n}^*)^{(-q/2)+\theta}), \tag{C.4}$$

where $0 \leq \theta = \theta(\xi) < 1$ is any exponent obtained from (K.6) with $\xi$ satisfying (4.3), and

$$E\left|\frac{\mathbf{f}(\hat{k}_{n,P_n}^S) - \mathbf{f}(k_{n,D_n}^*)}{(L_{n,D_n}(\hat{k}_{n,P_n}^S))^{1/2}}\right|^{2q} = o((P_n - 1)^{-q}). \tag{C.5}$$

Consequently, the claimed result is guaranteed by (3.6) [with $\xi$ satisfying (4.3)], (C.3)–(C.5) and Proposition 2. □

PROOF OF THEOREM 3. It suffices to show that (C.3)–(C.5) hold with $\hat{k}_{n,P_n}^S$ replaced by $\hat{k}_{n,P_n}$. Define $G_n(k) = N\exp\{\mathrm{IC}_{P_n}(k)\} - S_n^{(P_n)}(k)$ and $|U_{6,n}(k)| = |G_n(k) - G_n(k_{n,D_n}^*)|/NL_{n,D_n}(k)$. Then, by the same reasoning as in (C.2), $P(\hat{k}_{n,P_n} = k) \leq \sum_{i=1}^6 P(|U_{in}(k)| \geq (1/6)V_{n,D_n}(k))$. Moreover, Taylor's theorem and [14], (5.42), yield that for $q > 0$, all $1 \leq k \leq K_n$, and all sufficiently large $n$, $E|U_{6,n}(k)|^q \leq CP_n^{2q}K_n^{2q}N^{-2q}L_{n,D_n}^{-q}(k)$. These inequalities and the same argument used in the proof of Theorem C.1 give the desired results. □

## APPENDIX D: PROOFS OF THEOREMS 4 AND 6

PROOF OF THEOREM 4. First observe that for all sufficiently large $n$,

(D.1)
$$p_0 \leq K_n, \qquad k_n^* = k_{n,H_n}^* = p_0,$$
$$L_n(k_n^*) = \frac{p_0 \sigma^2}{N}, \qquad L_{n,H_n}(k_{n,H_n}^*) = \frac{(H_n - 1)p_0 \sigma^2}{N},$$

where $H_n = o(n)$ and $H_n > 1$. Define $I_{3,n}(k) = \{N(\mathbf{f}(k) - \mathbf{f}(p_0))^2\}/(p_0 \sigma^2)$. In view of (D.1) and by Hölder's inequality, we have, for all sufficiently large $n$,

(D.2)
$$(D_n - 1)\frac{E(\mathbf{f}(\hat{k}_{n,\delta_n}) - \mathbf{f}(k_{n,D_n}^*))^2}{L_{n,D_n}(k_{n,D_n}^*)}$$
$$\leq \sum_{k=1}^{p_0-1} (E|I_{3,n}(k)|^r)^{1/r} P^{(r-1)/r}(\hat{k}_{n,\delta_n} = k)$$
$$+ \sum_{k=p_0+1}^{K_n} (E|I_{3,n}(k)|^r)^{1/r} P^{(r-1)/r}(\hat{k}_{n,\delta_n} = k) \equiv (\mathrm{I}) + (\mathrm{II}),$$

where $r > 1$ and $D_n = D_{\mathrm{APE}_{\delta_n}}$. According to [14], Proposition 1 and Lemmas 1 and 2, (B.7) and (B.11),

(D.3)
$$E|I_{3,n}(k)|^r \leq \begin{cases} C, & 1 \leq k \leq p_0 - 1, \\ Ck^r, & p_0 + 1 \leq k \leq K_n, \end{cases}$$

as $n$ is sufficiently large. Armed with (D.3), Lemmas B.2–B.4, the conditions imposed on $\delta_n$, and the fact that for $1 \leq k \leq K_n$ and $k \neq p_0$, $V_{n,D_n}^{-1}(k) \leq C$, the proof of Corollary B.1 is modified to obtain that for any $s > 0$, $(\mathrm{I}) = O(n^{-s})$ and $(\mathrm{II}) = o((\log \delta_n^{-1})^{-s})$. Hence, (2.9) holds for $\hat{k}_{n,OS} = \hat{k}_{n,\delta_n}$. Similarly, it can be shown that $\hat{k}_{n,\delta_n}$ also satisfies (2.8). In view of Proposition 2 and (D.1), (2.3) [(2.4)] is achieved by $\hat{k}_{n,\delta_n}$. Modifying the proof of Theorem 3 and the above argument for $\mathrm{APE}_{\delta_n}(k)$, it can be shown that (2.8), (2.9) and (2.3) [(2.4)] can also be verified for $\hat{k}_{n,P_n}$, with $P_n$ satisfying the imposed constraints, $P_n \to \infty$ as $n \to \infty$ and $P_n = O(n^s)$ for some $0 < s < 1$. The details are omitted in order to save space. □

PROOF OF THEOREM 6. Unlike the previous theorems, Proposition 2 is not applied in the proof of Theorem 6 since the penalty term associated with $\hat{k}_n^{(\iota)}$, $2I_{\{\hat{k}_{n,P_n} \neq \hat{k}_{n^\iota,P_{n^\iota}}\}} + P_n I_{\{\hat{k}_{n,P_n} = \hat{k}_{n^\iota,P_{n^\iota}}\}}$, is random. In the following, we shall directly verify that

$$\text{(D.4)} \quad \limsup_{n\to\infty} \frac{q_n(\hat{k}_n^{(\iota)}) - \sigma^2}{L_n(k_n^*)} \leq 1.$$

First assume that condition (ii) of Theorem 6 holds. Let $0 < \xi < \min\{\delta_1^*/2, 1/2, (1+\delta_1^*) - \iota_1(2+\delta_1^*)\}$. Then there are $0 \leq \theta = \theta(\xi) < 1$ and $M = M(\xi) > 0$ such that (5.2) is satisfied. Define

$$B_{n,M^*} = A_{P_n,\theta,M}^C \cap \left\{ k : 1 \leq k \leq K_n, \frac{L_{n,P_n}(k) - L_{n,P_n}(k_{n,P_n}^*)}{L_{n,P_n}(k_{n,P_n}^*)} < M^* \right\},$$

where $M^*$ is some positive constant, and $I_{4,n}(k) = (\mathbf{f}(k) + \mathcal{S}(k))^2/L_n(k_n^*)$. Then

$$\frac{q_n(\hat{k}_n^{(\iota)}) - \sigma^2}{L_n(k_n^*)}$$

$$= E\{I_{4,n}(\hat{k}_n^{(\iota)})\}$$

$$\text{(D.5)} \quad \leq E\{I_{4,n}(\hat{k}_{n,2})\}$$

$$+ E\{I_{4,n}(\hat{k}_{n,P_n}) I_{\{\hat{k}_{n,P_n} = \hat{k}_{n^\iota,P_{n^\iota}}\}} (I_{\{\hat{k}_{n,P_n} \in B_{n,M^*}\}} + I_{\{\hat{k}_{n,P_n} \notin B_{n,M^*}\}})\}$$

$$\equiv (\text{I}) + (\text{II}).$$

Observe that for $r > 1$,

$$(\text{II}) \leq E\{I_{4,n}(\hat{k}_{n^\iota,P_{n^\iota}}) I_{\{\hat{k}_{n^\iota,P_{n^\iota}} \in B_{n,M^*}\}}\} + E\{I_{4,n}(\hat{k}_{n,P_n}) I_{\{\hat{k}_{n,P_n} \notin B_{n,M^*}\}}\}$$

$$= \sum_{\substack{k=1 \\ k \in B_{n,M^*}}}^{K_{n^\iota}} E\{I_{4,n}(k) I_{\{\hat{k}_{n^\iota,P_{n^\iota}} = k\}}\} + \sum_{\substack{k=1 \\ k \notin B_{n,M^*}}}^{K_n} E\{I_{4,n}(k) I_{\{\hat{k}_{n,P_n} = k\}}\}$$

$$\text{(D.6)} \quad \leq C \left\{ \sum_{\substack{k=1 \\ k \in B_{n,M^*}}}^{K_{n^\iota}} \frac{L_n(k)}{L_n(k_n^*)} P^{(r-1)/r}(\hat{k}_{n^\iota,P_{n^\iota}} = k) \right.$$

$$\left. + \sum_{\substack{k=1 \\ k \notin B_{n,M^*}}}^{K_n} \frac{L_n(k)}{L_n(k_n^*)} P^{(r-1)/r}(\hat{k}_{n,P_n} = k) \right\}$$

$$\equiv C\{(\text{III}) + (\text{IV})\},$$

where the second inequality follows from Hölder's inequality and the fact that for all $1 \leq k \leq K_n$, $E|\mathbf{f}(k) + \mathcal{S}(k)|^{2r} \leq CL_n^r(k)$, which is ensured by Wei ([28], Lemma 2), Ing and Wei ([14], Proposition 1 and Lemmas 1 and 2) and

(B.16). To deal with (III), by (5.3) and the definition of $B_{n,M^*}$, we have, for all sufficiently large $n$,

$$B_{n,M^*} \cap \{1, 2, \ldots, K_{n^\iota}\} \subseteq A_{P_{n^\iota}, \theta, M}.$$

Hence, (5.2) ensures that for all $k \in B_{n,M^*} \cap \{1, 2, \ldots, K_{n^\iota}\}$ and sufficiently large $n$,

(D.7)
$$V_{n^\iota, P_{n^\iota}}^{-1}(k) = \{L_{n^\iota, P_{n^\iota}}(k)\}\{L_{n^\iota, P_{n^\iota}}(k) - L_{n^\iota, P_{n^\iota}}(k_{n^\iota, P_{n^\iota}}^*)\}^{-1} \leq C(k_{n^\iota, P_{n^\iota}}^*)^\xi.$$

The definition of $B_{n,M^*}$ also yields for all $k \in B_{n,M^*}$ and $P_n \geq 2$, $L_n(k)/L_n(k_n^*) \leq (P_n - 1)\{L_{n,P_n}(k)/L_{n,P_n}(k_{n,P_n}^*)\} \leq CP_n$. According to this, (D.7) and the moment bounds for $|U_{in}(k)|, i = 1, \ldots, 6$, we have, for $q > 0$ and all sufficiently large $n$,

(D.8)
$$(\text{III}) \leq CP_n(k_{n^\iota, P_{n^\iota}}^*)^{\xi q} \Bigg\{ \sum_{\substack{k=1 \\ k \in B_{n,M^*}}}^{K_{n^\iota}} \frac{P_{n^\iota}^q(k^q + k_{n^\iota, P_{n^\iota}}^{*q})}{N_\iota^q} + \frac{k^{2q} + k_{n^\iota, P_{n^\iota}}^{*2q}}{N_\iota^{3q/2} L_{n^\iota, P_{n^\iota}}^q(k)} + \frac{1}{N_\iota^{q/2}} + \frac{k^{q/2} + k_{n^\iota, P_{n^\iota}}^{*q/2}}{N_\iota^q L_{n^\iota, P_{n^\iota}}^q(k)} + \frac{\|\mathbf{a}(k) - \mathbf{a}(k_{n^\iota, P_{n^\iota}}^*)\|_R^q}{N_\iota^{q/2} L_{n^\iota, P_{n^\iota}}^q(k)} + \frac{K_{n^\iota}^{2q} P_{n^\iota}^{2q}}{N_\iota^{2q} L_{n^\iota, P_{n^\iota}}^q(k)} \Bigg\}.$$

By taking $q$ on the right-hand side of (D.8) large enough and in view of the restrictions on $\xi$ (given at the beginning of this proof), $\iota$ and $P_n$,

(D.9) $\qquad (\text{III}) = o(1).$

Similarly, (5.2) and the definition of $B_{n,M^*}$ imply that for all sufficiently large $n$ and $k \in \{k : 1 \leq k \leq K_n \text{ and } k \notin B_{n,M^*}\}$, (D.7) is still valid if $n^\iota$ and $P_{n^\iota}$ are replaced by $n$ and $P_n$, respectively. This finding, the fact that for $P_n \geq 2$, $L_n(k)/L_n(k_n^*) \leq (P_n - 1)\{L_{n,P_n}(k)/L_{n,P_n}(k_{n,P_n}^*)\}$, and an argument similar to the one used to verify (D.9) give $(\text{IV}) = o(1)$, which, together with (D.5), (D.6), (D.9) and Theorem 2 of [14], yields (D.4).

Next, assume that condition (i) holds. By the moment bounds for $|U_{in}(k)|$, $i = 1, \ldots, 6$, and similar reasoning to that used in the proof of Theorem 4, we have

(D.10) $\qquad \lim_{n \to \infty} P(\hat{k}_{n,P_n} \neq p_0) = 0,$

and for any $q > 0$,

(D.11) $\qquad E|I_{4,n}(\hat{k}_{n,2})|^q = O(1).$

Since

$$\frac{q_n(\hat{k}_n^{(\iota)}) - \sigma^2}{L_n(k_n^*)} \leq E\{I_{4,n}(\hat{k}_{n,P_n})\}$$

(D.12)
$$+ E\{I_{4,n}(\hat{k}_{n,2})(I_{\{\hat{k}_{n,P_n} \neq p_0\}} + I_{\{\hat{k}_{n^\iota,P_{n^\iota}} \neq p_0\}})\},$$

(D.4) follows from (D.10)–(D.12), Hölder's inequality and $\limsup_{n\to\infty} E(I_{4,n}(\hat{k}_{n,P_n})) \leq 1$ (which is ensured by Theorem 4). This completes the proof of the theorem. $\square$

**Acknowledgments.** I would like to thank a Co-Editor, an Associate Editor, and two anonymous referees for their valuable and constructive comments which helped to substantially improve the readability and content of the paper. I am also grateful to Atsushi Inoue and Lutz Kilian for kindly correcting some errors in my GAUSS code for implementing the simulations in a previous version of this paper.

Institute of Statistical Science
Academia Sinica 128
Academia Road, Section 2
Taipei 115
Taiwan
Republic of China
E-mail: cking@stat.sinica.edu.tw